\documentclass[12pt,reqno]{amsart} 

\usepackage[all]{xy}
\usepackage{supertabular}

\def\theoremnumberstyle{plain} 
\def\proofname{Proof}
\usepackage{tom}

\newcommand{\punkte}{\lpatt (0.067 0.1) \setgray 0.5 }
\newcommand{\strich}{\lpatt () \setgray 0 }
\newcommand{\aus}[2]{#1}
\newcommand{\muell}[1]{}

\numberwithin{equation}{section} 
\allowdisplaybreaks[3]   

\SelectTips{cm}{12}  
\CompileMatrices

\renewcommand{\arraystretch}{1.4}

\begin{document}

   \parindent0cm

   \title[The $\gamma_\alpha^*$-Invariant]{A new Invariant for Plane Curve Singularities}
   \author{Thomas Keilen}
   \address{University of Warwick\\
     Mathematics Institute\\
     Coventry CV4 7AL\\
     United Kingdom
     }
   \email{keilen@mathematik.uni-kl.de}
   \urladdr{http://www.mathematik.uni-kl.de/\textasciitilde keilen}
   \author{Christoph Lossen}
   \address{Universit\"at Kaiserslautern\\
     Fachbereich Mathematik\\
     Erwin-Schr\"odinger-Stra\ss e\\
     D -- 67663 Kaiserslautern\\
     Germany
     }
   \email{lossen@mathematik.uni-kl.de}
   \urladdr{http://www.mathematik.uni-kl.de/\textasciitilde lossen}
   \thanks{The first author was supported by the European
     mathematical network EAGER. The second author was supported by
     DFG grant no. Lo 864/1.}

   \subjclass{Primary 14H20, 14H15, 14H10}

   \date{September, 2003.}

   \keywords{Algebraic geometry, singularity theory}
     
   \begin{abstract}
     In \cite{GLS01} the authors gave a general sufficient
     numerical condition for the T-smoothness (smoothness and expected
     dimension) of equisingular
     families of plane curves. This condition involves a new invariant
     $\gamma^*$ for
     plane curve singularities, and it is conjectured to be
     asymptotically proper. In \cite{Kei04}, similar sufficient
     numerical conditions are obtained for the T-smoothness of
     equisingular families on various classes surfaces. These
     conditions involve a series of  invariants
     $\gamma_\alpha^*$, $0\leq \alpha\leq 1$, with
     $\gamma_1^*=\gamma^*$.
     In the present paper we compute (respectively give bounds for)
     these invariants for semiquasihomogeneous singularities.
   \end{abstract}

   \maketitle

   When studying numerical conditions for the T-smoothness of
   equisingular families of curves, new  invariants of plane curve
   singularities $V(f)\subset(\C^2,0)$ turn up. These invariants are
   defined as the maximum of a function depending on the
   codimension of complete intersection ideals  containing the
   Tjurina ideal, respectively the equisingularity ideal, of $f$, and on the
   intersection multiplicity of $f$
   with elements of the complete intersection ideals. In Section \ref{sec:invariants}
   we will define these invariants, and we will calculate
   them for several classes of singularities, the main results being
   Proposition \ref{prop:gamma-simple}, Proposition
   \ref{prop:gamma-ordinary} and Proposition \ref{prop:gamma}. It is
   the upper bound in Lemma \ref{lem:gamma} which ensures that the
   conditions for T-smoothness with these new conditions (see
   \cite{GLS00}, \cite{GLS01}, \cite{Kei04}) improve
   than the previously known ones (see \cite{GLS97}). In
   the remaining sections we introduce some notation and we gather
   some necessary, though mainly well-known technical results used in the proofs of Section
   \ref{sec:invariants}. 

   We should like to point out that the definition of the invariant
   $\gamma_1^*$ below is a modification of the invariant ``$\gamma^*$''
   defined in \cite{GLS01}, and it is always bound from above by the
   latter. Moreover, the latter can be replaced by it in the
   conditions of \cite{GLS01} Proposition 2.2.

   \begin{varthm-roman-break}[Notation]
     Throughout this paper $R=\C\{x,y\}$ will be the ring of
     convergent power series in the variables $x$ and $y$, and
     $\m=\langle x,y\rangle\lhd R$ will be its maximal ideal.
   \end{varthm-roman-break}

   \tableofcontents

   \section{The $\gamma_\alpha^*$-Invariants}\label{sec:invariants}

   For the definition of the $\gamma_\alpha^*$-invariants the Tjurina
   ideal, respectively the equisingularity ideal in the sense
   of \cite{Wah74}, play an essential role. For the convenience of the
   reader we recall their definitions. 

   \begin{definition}\label{def:tjurina}
     Let $f\in \m$ be a reduced power series. The \emph{Tjurina ideal} of $f$
     is defined as 
     \begin{displaymath}
       I^{ea}(f)=\left\langle\frac{\partial f}{\partial
           x},\frac{\partial f}{\partial y},f \right\rangle,
     \end{displaymath}
     and the \emph{equisingularity ideal} of $f$ is defined as
     \begin{displaymath}
       I^{es}(f)=\big\{g\in R\;\big|\;f+\varepsilon g \mbox{ is equisingular
         over } \C[\varepsilon]/(\varepsilon^2)\big\}\supseteq I^{ea}(f).
     \end{displaymath}
     Their codimensions
     \begin{displaymath}
       \tau(f)=\dim_\C R/I^{ea}(f),
     \end{displaymath}
     respectively
     \begin{displaymath}
       \tau^{es}(f)=\dim_\C R/I^{es}(f),
     \end{displaymath}
     are analytical, respectively topological, invariants of
     the singularity type defined by $f$. Note that $\tau^{es}(f)$ is
     the codimension of the $\mu$-constant stratum in the equisingular
     deformation of the plane curve singularity defined by $f$. It can
     be computed in terms of multiplicities of the strict transform of
     $f$ at essential infinitely near points in the resolution tree of
     $\big(V(f),0\big)$ (cf.\ \cite{Shu91a}).
   \end{definition}

   \begin{definition}
     Let  $f\in\m$ be a reduced power
     series, and let $0\leq \alpha\leq 1$ be a rational number.
     \\
     If $I$ is a  zero-dimensional ideal in $R$ with
     $I^{ea}(f)\subseteq I\subseteq \m$ and $g\in I$, we define 
     \begin{displaymath}
       \lambda_\alpha(f;I,g):=\frac{\big(\alpha\cdot
         i(f,g)+(1-\alpha)\cdot \dim_\C(R/I)\big)^2}{i(f,g)-\dim_\C(R/I)},
     \end{displaymath}
     and
     \begin{displaymath}
       \gamma_\alpha(f;I):=\max\left\{(1+\alpha)^2\cdot\dim_\C(R/I),\;\lambda_\alpha(f;I,g)\;
         \big|\; g\in I, i(f,g)\leq 2\cdot\dim_\C(R/I)\right\},
     \end{displaymath}
     where $i(f,g)$ denotes the intersection multiplicity of $f$ and
     $g$. Note that, by Lemma \ref{lem:schnittzahl}, $i(f,g)>\dim_\C(R/I)$
     for all $g\in I$. Thus $\gamma_\alpha(f;I)$ is a well-defined positive rational
     number.
     \\
     We then set
     \begin{displaymath}
       \;\;\;\;\;\gamma_\alpha^{ea}(f):=\max\left\{0,\;\gamma_\alpha(f;I)\;\big|\;I\supseteq
         I^{ea}(f) \mbox{ is a complete intersection ideal}\right\}
     \end{displaymath}
     and
     \begin{displaymath}
       \;\;\;\;\;\gamma_\alpha^{es}(f):=\max\left\{0,\;\gamma_\alpha(f;I)\;\big|\;I\supseteq
         I^{es}(f) \mbox{ is a complete intersection ideal}\right\}
     \end{displaymath}
     Note that if $f\in\m\setminus\m^2$, then
     $I^{ea}(f)=I^{es}(f)=R$ and there is no zero-dimensional complete
     intersection ideal containing any of those two, hence
     $\gamma_\alpha^{ea}(f)=\gamma_\alpha^{es}(f)=0$. 
   \end{definition}

   \begin{lemma}\label{lem:schnittzahl}
     Let $f\in\m^2$ be reduced, and let
     $I$ be an ideal such that $I^{ea}(f)\subseteq I\subseteq \m$. 
     \\
     Then, for any $g\in I$,
     we have
     \begin{displaymath}
       \dim_\C(R/I)<\dim_\C\big(R/\langle f,g\rangle\big)=i(f,g).
     \end{displaymath}
   \end{lemma}
   \begin{proof}
     Cf.~\cite{Shu97} Lemma 4.1; the idea is mainly to show that not both
     derivatives of $f$ can belong to $\langle f,g\rangle$.
   \end{proof}

   Up to embedded isomorphism the Tjurina ideal 
   only depends on the analytical type of the
   singularity. More precisely, if $f\in R$ any power series, $u\in R$ a unit and
   $\phi:R\rightarrow R$ an isomorphism, then
   $I^{ea}(u\cdot f\circ\phi)=\{g\circ\phi\;|\;g\in I^{ea}(f)\}$.
   Thus the following definition makes sense.

   \begin{definition}
     Let $\ks$ be an analytical, respectively topological, singularity
     type, and let $f\in R$ be a representative of $\ks$. We
     then define 
     \begin{displaymath}
       \gamma_\alpha^{ea}(\ks):=\gamma_\alpha^{ea}(f),
     \end{displaymath}
     respectively 
     \begin{displaymath}
       \gamma_\alpha^{es}(\ks):=\max\{\gamma_\alpha^{es}(g)\;|\;g\mbox{ is a representative of }\ks\}.
     \end{displaymath}
   \end{definition}

   Since $i(f,g)>\dim_\C(R/I)$ in the above situation, we deduce the
   following lemma.

   \begin{lemma}\label{lem:gam}
     Let $f\in\m^2$ be reduced, 
     $I^{ea}(f)\subseteq I\subseteq \m$ be a zero-dimensional
     ideal, and $0\leq \alpha<\beta\leq 1$, then
     $\gamma_\alpha(f;I)< \gamma_\beta(f;I)$.

     In particular, for any analytical, respectively topological,
     singularity type 
     \begin{displaymath}
       \gamma_\alpha^{ea}(\ks)<\gamma_{\beta}^{ea}(\ks)
       \;\;\;\;\;\mbox{ respectively }\;\;\;\;\;
       \gamma_\alpha^{es}(\ks)<\gamma_{\beta}^{es}(\ks).
     \end{displaymath}     
   \end{lemma}

   For reasons of comparison let us also recall the definition of 
   $\tau_{ci}^{ea}$, $\tau_{ci}^{es}$, $\kappa$ and $\delta$.

   \begin{definition}
     For $f\in R$ we define
     \begin{displaymath}
       \tau_{ci}^{ea}(f):=\max\{0,\dim_\C(R/I)\;|\;I\supseteq I^{ea}(f)
       \mbox{ a complete intersection}\},
     \end{displaymath}
     and 
     \begin{displaymath}
       \tau_{ci}^{es}(f):=\max\{0,\dim_\C(R/I)\;|\;I\supseteq I^{es}(f)
       \mbox{ a complete intersection}\}.
     \end{displaymath}
     Again, for analytically equivalent
     singularities the values coincide, so that for an analytical
     singularity type $\ks$, choosing some
     representative $f\in R$, we may define
     \begin{displaymath}
       \tau_{ci}^{ea}(\ks):=\tau_{ci}(f).
     \end{displaymath}
     For a topological singularity type we set
     \begin{displaymath}
       \tau_{ci}^{es}(\ks):=\max\{\tau_{ci}^{es}(g)\;|\;g\mbox{ a
         representative of }\ks\}.
     \end{displaymath}
     Note that obviously
     \begin{displaymath}
       \tau_{ci}^{ea}(\ks)\leq\tau(\ks)\;\;\;\mbox{ and }\;\;\;
       \tau_{ci}^{es}(\ks)\leq\tau^{es}(\ks),
     \end{displaymath}
     where $\tau(\ks)$ is the Tjurina number of $\ks$ and
     $\tau^{es}(\ks)$ is  as
     defined in Definition \ref{def:tjurina}.
   \end{definition}

   \begin{definition}
     For $f\in R$ and $\ko=R/\langle f\rangle$, we define the $\delta$-invariant
     \begin{displaymath}
       \delta(f)=\dim_\C \widetilde{\ko}/\ko
     \end{displaymath}
     where $\ko\; \subset\; \widetilde{\ko}$ is the normalisation of
     $\ko$, and the $\kappa$-invariant
     \begin{displaymath}
       \kappa(f)=i\left(f,\alpha\cdot\frac{\partial f}{\partial x}+\beta\cdot\frac{\partial f}{\partial x}\right),
     \end{displaymath}
     where $(\alpha:\beta)\in\PC^1$ is generic.
       
     $\delta$ and $\kappa$ are topological (thus also analytical) invariants
     of the singularity defined by $f$ so that for the topological,
     respectively analytical, singularity type $\ks$ given by $f$ we
     can set
     \begin{displaymath}
       \delta(\ks)=\delta(f)\;\;\;\mbox{ and }\;\;\;\kappa(\ks)=\kappa(f).
     \end{displaymath}     
   \end{definition}

   \medskip
   \begin{center}
     \framebox[13cm]{
       \begin{minipage}{12cm}
         \medskip
         Throughout this article we will sometimes treat
         topological and analytical singularities at the same time.
         Whenever we do so, we will write $I^*(f)$ for $I^{ea}(f)$
         respectively for $I^{ea}(f)$, and analogously we will use the
         notation $\gamma_\alpha^*$, $\tau_{ci}^*$ and $\tau^*$.
         \medskip
       \end{minipage}
       }
   \end{center}
   \bigskip

   The following lemma is again obvious from the definition of
   $\gamma_\alpha(f;I)$, once we take into account that $\kappa(f)=i(f,g)$ for
   a generic element $g\in I^{ea}(f)$ of $f$ and that for a fixed value of
   $d=\dim_\C(R/I)$ the function $i\mapsto\frac{(\alpha
     i+(1-\alpha)\cdot d)^2}{i-d}$ takes its maximum on $[d+1,2d]$ for
   the minimal possible value $i=d+1$. 

   \begin{lemma}\label{lem:gamma}
     Let $f\in\m^2$ be reduced, and let
     $I$ be an ideal in $R$ such that $I^{ea}(f)\subseteq I\subseteq \m$.

     Then
     \begin{displaymath}
       (1+\alpha)^2\cdot \dim_\C(R/I)\leq \gamma_\alpha(f;I)\leq \big(\dim_\C(R/I)+\alpha\big)^2.
     \end{displaymath}
     Moreover, if $\kappa(f)\leq 2\cdot\dim_\C(R/I)$, then
     \begin{displaymath}
       \gamma_\alpha(f;I)\geq
       \frac{\big(\alpha\cdot
         \kappa(f)+(1-\alpha)\cdot\dim_\C(R/I)\big)^2}{\kappa(f)-\dim_\C(R/I)}.
     \end{displaymath}
     In particular, for any analytical, respectively topological,
     singularity type $\ks$
     \begin{displaymath}
       (1+\alpha)^2\cdot\tau_{ci}^*(\ks)\leq \gamma_\alpha^*(\ks)\leq\big(\tau_{ci}^*(\ks)+\alpha\big)^2,
     \end{displaymath}
     and if $\kappa(\ks)\leq 2\cdot\tau^*_{ci}(\ks)$, then
     \begin{displaymath}
       \gamma_\alpha^*(\ks)\geq
       \frac{\big(\alpha\cdot
         \kappa(\ks)+(1-\alpha)\cdot\tau_{ci}^*(\ks)\big)^2}{\kappa(\ks)-\tau_{ci}^*(\ks)}.
     \end{displaymath}
   \end{lemma}

   In order to make the conditions for T-smoothness in \cite{Kei04} as sharp as
   possible, it is useful to know under which circumstances the term
   $(1+\alpha)^2\cdot \dim_\C(R/I)$ involved in the definition of
   $\gamma_\alpha^*(f)$ is actually exceeded. 

   \begin{lemma}
     If $\ks$ is a topological or analytical singularity type such
     that $\kappa(\ks)<2\cdot\tau_{ci}^*(\ks)$, then
     \begin{displaymath}
       (1+\alpha)^2\cdot\tau^*_{ci}(\ks)<\gamma_\alpha^*(\ks).
     \end{displaymath}
     This is in particular the case, if $\ks\not=A_1$ and $\tau_{ci}^*(\ks)=\tau^*(\ks)$, i.\ e.\ if
     the Tjurina ideal, respectively the equisingularity ideal, of some
     representative is a complete intersection.

   \end{lemma}
   \begin{proof}
     Lemma
     \ref{lem:gamma} gives
     \begin{displaymath}
       \gamma_\alpha^*(\ks)\geq \frac{\big(\alpha\cdot \kappa(\ks) +
         (1-\alpha) \cdot\tau_{ci}^*(\ks)\big)^2}{\kappa(\ks)-\tau_{ci}^*(\ks)}.
     \end{displaymath}
     If we consider the right-hand side as a function in $\kappa(\ks)$,
     it is strictly decreasing  on the interval
     $\left[0,2\cdot\tau^*_{ci}(\ks)\right]$ and takes its minimum thus
     at $2\cdot\tau^*_{ci}(\ks)$. By the assumption on $\kappa(\ks)$ we,
     therefore, get
     \begin{displaymath}
       \gamma_\alpha^*(\ks)>(1+\alpha)^2\cdot\tau_{ci}^*(\ks).
     \end{displaymath}
     Suppose now that $\tau^*_{ci}(\ks)=\tau^*(\ks)$ and $\ks\not=A_1$. By Lemma 
     \ref{lem:delta} we know
     $\delta(\ks)<\tau^{es}(\ks)\leq\tau(\ks)$. On the other hand 
     we have $\kappa(\ks)\leq
     2\cdot\delta(\ks)$ (see  \cite{GLS05}). 
     Therefore, $\kappa(\ks)<2\cdot\tau_{ci}^*(\ks)$. 
   \end{proof}

   \begin{lemma}\label{lem:delta}
     If $\ks\not=A_1$ is any analytical or topological singularity type, then
     $\delta(\ks)<\tau^{es}(\ks)$.
   \end{lemma}
   \begin{proof}
     If $(C,z)$ is a representative of $\ks$ and if $\kt^*(C,z)$ is
     the essential subtree of the complete embedded resolution tree of
     $(C,z)$, then 
     \begin{displaymath}
       \delta(\ks)=\sum_{p\in\kt^*(C,z)}\frac{\mult_p(C)\cdot(\mult_p(C)-1)}{2}
     \end{displaymath}
     and
     \begin{displaymath}
       \tau^{es}(\ks)=\sum_{p\in\kt^*(C,z)}\frac{\mult_p(C)\cdot(\mult_p(C)+1)}{2}
       -\#\mbox{ free points in }\kt^*(C,z)-1,
     \end{displaymath}
     where $\mult_p(C)$ denotes the multiplicity of the strict
     transform of $C$ at $p$ (see \cite{GLS05}). 
     Setting $\varepsilon_p=0$ if $p$ is
     satellite, $\varepsilon_p=1$ if $p\not=z$ is free, and $\varepsilon_z=2$,
     then $\mult_p(C)\geq\varepsilon_p$ and therefore
     \begin{displaymath}
       \tau^{es}(\ks)=\delta(\ks)+\sum_{p\in\kt^*(C,z)}\big(\mult_p(C)-\varepsilon_p\big)\geq \delta(\ks). 
     \end{displaymath}
     Moreover, we have equality if and only if $\mult_z(C)=2$,
     $\mult_p(C)=1$ for all $p\not=z$ and there is no satellite
     point, but this implies that $\ks=A_1$.
   \end{proof}

   For some classes of singularities  we can calculate the
   $\gamma_\alpha^*$-invariant concretely, and for some
   others we can at least give an upper bound, which in general is much
   better than the one derived from Lemma \ref{lem:gamma}. We restrict
   our attention to singularities
   having a convenient semi-quasihomogeneous 
   representative $f\in R$ (see Definition
   \ref{def:sqh}). Throughout the following proofs we will frequently
   make use of monomial orderings, see Section~\ref{sec:monomialorderings}. 

   \begin{proposition}[(Simple Singularities)]\label{prop:gamma-simple}
     Let $\alpha$ be a rational number with $0\leq \alpha\leq 1$.
     Then we obtain the following values for
     $\gamma_\alpha^{es}(\ks)=\gamma_\alpha^{ea}(\ks)$, where $\ks$ is
     a simple singularity type.
     \begin{displaymath}
       \begin{array}{|cc|c|}
         \hline
         \multicolumn{2}{|c|}{\ks}  & \gamma_\alpha^{ea}(\ks)= \gamma_\alpha^{es}(\ks)
         \\\hline\hline
         A_k, & k\geq 1 & (k+\alpha)^2 \\\hline
         D_k, & 4\leq k\leq 4+\sqrt{2}\cdot(2+\alpha) 
         & \frac{(k+2\alpha)^2}{2} \\\hline
         D_k, & k\geq 4+\sqrt{2}\cdot(2+\alpha) & (k-2+\alpha)^2\\\hline
         E_k, & k=6,7,8 & \frac{(k+2\alpha)^2}{2}  \\\hline
       \end{array}
     \end{displaymath}
   \end{proposition}
   \begin{proof}
     Let $\ks_k$ be one of the simple singularity types $A_k$,
     $D_k$ or $E_k$, and let $f\in R$ be a representative of
     $\ks_k$. 
     Note that the Tjurina ideal $I^{ea}(f)$ and the equisingularity ideal
     $I^{es}(f)$ coincide, and hence so do the $\gamma_\alpha^*$-invariants, i.\
     e.\ 
     \begin{displaymath}
       \gamma_\alpha^{ea}(\ks_k)=\gamma_\alpha^{es}(\ks_k).
     \end{displaymath}

     Moreover, in the considered cases the Tjurina ideal
     is indeed a complete intersection ideal with
     $\dim_\C\big(R/I^{ea}(f)\big)=k$, so that in particular the given
     values are upper bounds for $(1+\alpha)^2\cdot\dim_\C(R/I)$ for
     any complete intersection ideal $I$ containing the Tjurina ideal.
     By Lemma \ref{lem:gamma}
     we know
     \begin{displaymath}
       \frac{(\alpha\cdot
         \kappa(\ks_k)+(1-\alpha)\cdot k)^2}{\kappa(\ks_k)-k}\leq 
       \gamma_\alpha(\ks_k)\leq  (k+\alpha)^2.
     \end{displaymath}
     Note that $\kappa(A_k)=k+1$, $\kappa(D_k)=k+2$ and
     $\kappa(E_k)=k+2$, which in particular gives the result for
     $\ks_k=A_k$. Moreover, it shows that for $\ks_k=D_k$ or
     $\ks_k=E_k$ we have
     \begin{displaymath}
       \gamma_\alpha(\ks_k)\geq \frac{(k+2\alpha)^2}{2}.
     \end{displaymath}

     If we fix a complete intersection ideal $I$ with
     $I^{ea}(f)\subseteq I$, then
     \begin{displaymath}
       \lambda_\alpha(f;I,g)=\frac{\big(\alpha\cdot i(f,g)+(1-\alpha)\cdot
         \dim_\C(R/I)\big)^2}{i(f,g)-\dim_\C(R/I)}, 
     \end{displaymath}
     with $g\in I$ such that $i(f,g)\leq 2\cdot\dim_\C(R/I)$, 
     considered as a function in $i(f,g)$ is maximal,
     when $i(f,g)$ is minimal. If $i(f,g)-\dim_\C(R/I)\geq 2$, then
     \begin{displaymath}
       \lambda_\alpha(f;I,g)\leq \frac{(k+2\alpha)^2}{2}.
     \end{displaymath}
     It therefore remains to consider the
     case where 
     \begin{equation}\label{eq:simple:1}
       i(f,g)-\dim_\C(R/I)=1
     \end{equation}
     for some $I$ and some $g\in
     I$, and to maximise the possible $\dim_\C(R/I)$.

     We claim that for $\ks_k=D_k$ with $f=x^2y-y^{k-1}$ as representative,
     $\dim_\C(R/I)\leq k-2$, and thus
     $I=\langle x,y^{k-2}\rangle$ and $g=x$ are suitable with
     \begin{displaymath}
       \lambda_\alpha(f;I,x)=(k-2+\alpha)^2,
     \end{displaymath}
     which is greater than $\frac{(k+2\alpha)^2}{2}$ if and only if
     $k\geq 4+\sqrt{2}\cdot (2+\alpha)$. Suppose, therefore,
     $\dim_\C(R/I)=k-1$. Then  $y^{k-1},x^3\in I^{ea}(f)=\langle
     xy,x^2-{\scriptsize(k-1)}\cdot y^{k-2}\rangle\subset I$, 
     the leading ideal $L_{<_{ls}}\big(I^{ea}(f)\big)=\langle
     x^3,xy,y^{k-2}\rangle\subset L_{<_{ls}}(I)$, and since by Proposition~
     \ref{prop:ordering} $\dim_\C(R/I)=\dim_\C\big(R/L_{<_{ls}}(I)\big)$,
     either $L_{<_{ls}}(I)=\langle x^3,xy,y^{k-3}\rangle$ or
     $L_{<_{ls}}(I)=\langle x^2,xy,y^{k-2}\rangle$. In the first case
     there is a power series $g\in I$ such that $g\equiv y^{k-3}+ax+bx^2\;(\mod
     I)$, and hence $I\ni yg\equiv y^{k-2}\;(\mod
     I)$, i.\ e.\ $y^{k-2}\in I$. But then $x^2\in
     I$ and $x^2\in L_{<_{ls}}(I)$, in contradiction to the assumption. In the second case,
     similarly, there is a $g\in I$ such that $g\equiv x^2\;(\mod
     I)$, and hence $x^2\in I$ which in turn implies that $y^{k-2}\in
     I$. Thus $I=\langle x^2,xy,y^{k-2}\rangle$, and $\dim_\C(I/\m
     I)=3$ which by Remark \ref{lem:nak} contradicts the
     fact that $I$ is a complete intersection.

     \aus{
       If $\ks_k=E_6$, then $f=x^3-y^4$ is a representative and
       $I^{ea}(f)=\langle x^2,y^3\rangle$. Suppose that
       $\dim_\C(R/I)=k-1=5$, then $L_{<_{ds}}(I)=\langle
       x^2,y^3,xy^2\rangle$ and $H^0_{R/I}=H^0_{R/L_{<_{ds}}(I)}$,
       in contradiction to Lemma \ref{lem:ci}, since
       $H^0_{R/L_{<_{ds}}(I)}(2)=2$ and $H^0_{R/L_{<_{ds}}(I)}(3)=0$.
       Thus $\dim_\C(R/I)\leq 4$ and $\lambda_\alpha(f;I,g)\leq
       (4+\alpha)^2\leq \frac{(6+2\alpha)^2}{2}$.
       
       If $\ks_k=E_7$, then $f=x^3-xy^3$ is a representative and 
       $I^{ea}(f)=\langle 3x^2-y^3,xy^2\rangle\ni x^3,y^5$. If
       $\dim_\C(R/I)\leq 4$, then $\lambda_\alpha(f;I,g)\leq
       (4+\alpha)^2\leq \frac{(7+2\alpha)^2}{2}$, and we are done. It
       thus remains to exclude the cases where
       $\dim_\C(R/I)\in\{5,6\}$.  For this we note first that if there
       is a $g\in I$ such that $L_{<_{ls}}(g)=y^2$, then 
       \begin{equation}\label{eq:simple:2}
         g\equiv
         y^2+ax+bx^2+cxy+dx^2y\;(\mod I),
       \end{equation}
       and therefore $y^2g\equiv
       y^4\;(\mod I)$, which implies $y^4\in I$ and hence $x^2y\in
       I$. Analogously, if there is a $g\in I$ such that
       $L_{<_{ls}}(g)=x^2y$, then $g\equiv x^2y\;(\mod I)$ and
       again $x^2y,y^4\in I$. Suppose now that $\dim_\C(R/I)=6$, then
       $L_{<_{ls}}(I)=\langle y^2,x^3\rangle$ or  
       $L_{<_{ls}}(I)=\langle y^3,xy^2,x^2y,x^3\rangle$. In both
       cases we thus have $x^2y,y^4\in I$. However, in the first case
       then $x^2y\in L_{<_{ls}}(I)$, in contradiction to the
       assumption. While in the second case we find $I=\langle 
       xy^2,x^2y,3x^2-y^3\rangle$, and $\dim_\C(I/\m
       I)=3$ contradicts the fact that $I$ is a complete intersection by
       Lemma \ref{lem:nak}.
       Suppose, therefore, that
       $\dim_\C(R/I)=5$. Then $L_{<_{ls}}(I)=\langle
       y^2,x^2y,x^3\rangle$, or $L_{<_{ls}}(I)=\langle
       y^3,xy^2,x^2\rangle$, or $L_{<_{ls}}(I)=\langle
       y^3,xy,x^3\rangle$. In the first case, we know already that
       $y^4,x^2y\in I$. Looking once more on \eqref{eq:simple:2} we
       consider the cases $a=0$ and $a\not=0$. If $a=0$, then $yg\equiv
       y^3\;(\mod I)$, and thus $y^3\in I$, which in turn implies
       $x^2\in I$. Similarly, if $a\not=0$, then $xg\equiv ax^2\;(\mod
       I)$ implies $x^2\in I$. But then also $x^2\in L_{<_{ls}}(I)$, in
       contradiction to the assumption. In the second case there is a
       $g\in I$ such that $g\equiv x^2+ax^2y\;(\mod I)$, and thus
       $yg\equiv x^2y\in I$. But then also $x^2\in I$ and $y^3\in I$, so
       that $I=\langle y^3,xy^2,x^2\rangle$. However, $\dim_\C(I/\m
       I)=3$ contradicts again the fact that $I$ is a complete intersection.
       Finally in the third case
       there is a $g\in I$ with $g\equiv xy+ax^2+bx^2y\;(\mod I)$, and thus
       $xg\equiv x^2y\;(\mod I)$ implies $x^2y\in I$ and then $xy+ax^2\in
       I$. Therefore, $I=\langle xy+ax^2,3x^2-y^3\rangle$, and for
       for $h\in I$ and for generic $b,c\in \C$  we have $i(f,h)\geq
       i(x,h)+i\big(x^2-y^3,b\cdot(xy+ax^2)+c\cdot(3x^2-y^3)\big)\geq 3+5=8$, in
       contradiction to \eqref{eq:simple:1}. 
       
       Finally, if $\ks_k=E_8$ with representative $f=x^3-y^5$ and
       $I^{ea}(f)=\langle x^2,y^4\rangle$, we get for
       $\dim_\C(R/I)\leq 5$ that $\lambda_\alpha(f;I,g)\leq
       (5+\alpha)^2\leq \frac{(8+2\alpha)^2}{2}$. It therefore remains
       to exclude the cases $\dim_\C(R/I)\in \{6,7\}$.
       If $\dim_\C(R/I)=7$ then $L_{<_{ds}}(I)=\langle
       x^2,y^4,xy^3\rangle$. But then $H^0_{R/L_{<_{ds}}(I)}(3)=2$
       and $H^0_{R/L_{<_{ds}}(I)}(4)=0$ are in contradiction to  Lemma
       \ref{lem:ci}. And if $\dim_\C(R/I)= 6$, then
       $L_{<_{ls}}(I)=\langle y^3,x^2\rangle$ or
       $L_{<_{ls}}(I)=\langle y^4,xy^2,x^2\rangle$.
       In the first case there is some $g\in I$ such that $g\equiv
       y^3+ax+bxy+cxy^2+dxy^3 \;(\mod I)$, and thus $xg\equiv
       xy^3\;(\mod I)$ and $xy^3\in I$. But then $yg\equiv
       axy+bxy^2\;(\mod I)$ and hence $axy+bxy^2\in I$. Since neither $xy\in L_{<_{ls}}(I)$ nor
       $xy^2\in L_{<_{ls}}(I)$, we must have $a=0=b$. Therefore,
       $g\equiv y^3+cxy^2\;(\mod I)$ and
       $I=\langle x^2,y^3+cxy^2\rangle$, which for $h\in I$ and
       $a,b\in\C$ generic  gives
       $i(f,g)\geq i\big(x^3-y^4,ax^2+b\cdot(y^3+cxy^2)\big)\geq 8$, in
       contradiction to \eqref{eq:simple:1}. In the second case, there
       is $g\in I$ such that $g\equiv xy^2+axy^3\;(\mod I)$, therefore
       $yg\equiv xy^3\;(\mod I)$ and $xy^3\in I$. But then $xy^2\in I$
       and $I=\langle y^4,xy^2,x^2\rangle$. This, however, is not a
       complete intersection, since $\dim_\C(I/\m I)=3$, in
       contradiction to the assumption.

       This finishes the proof.
       }{The cases of the exceptional singularities $E_6$, $E_7$ and
       $E_8$ are treated similarly.}
   \end{proof}

   \begin{proposition}[(Ordinary Multiple Points)]\label{prop:gamma-ordinary}
     Let $\alpha$ be a rational number with $0\leq \alpha\leq 1$, and
     let $M_k$ denote the topological singularity type of an ordinary
     $k$-fold point with $k\geq 3$. Then 
     \begin{displaymath}
       \gamma_\alpha^{es}(M_k)=2\cdot (k-1+\alpha)^2.
     \end{displaymath}
     In particular
     \begin{displaymath}
       \gamma_\alpha^{es}(M_k)>(1+\alpha)^2\cdot\tau^{es}_{ci}(M_k).
     \end{displaymath}
   \end{proposition}
   \begin{proof}
     Note that for any representative $f$ of $M_k$ we have
     $$I^{es}(f)=I^{ea}(f)+\m^k=\left\langle  
       \frac{\partial f_k}{\partial x},\frac{\partial f_k}{\partial y}
     \right\rangle+\m^k,$$ 
     where $f_k$ is the homogeneous part of
     degree $k$ of $f$, so that we may assume $f$ to be
     homogeneous of degree $k$. 
 
     If $I$ is a complete intersection ideal with
     $\m^k\subset I^{es}(f)\subseteq I$, then by Lemma \ref{lem:deg-fp}
     \begin{displaymath}
       \dim_\C(R/I)\leq \big(k-\mult(I)+1\big)\cdot\mult(I).
     \end{displaymath}
     We note moreover that for any $g\in I$ 
     \begin{displaymath}
       i(f,g)\geq\mult(f)\cdot\mult(g)\geq k\cdot\mult(I),
     \end{displaymath}
     and that for a fixed $I$ we may attain an upper bound for $\lambda_\alpha(f;I,g)$
     by replacing $i(f,g)$ by a lower bound for $i(f,g)$.

     Hence, if $\mult(I)\geq 2$, we have
     \begin{equation}\label{eq:multpoint:1}
       \lambda_\alpha(f;I,g)\leq 
       \frac{\big(k-(1-\alpha)\cdot(\mult(I)-1)\big)^2\cdot\mult(I)^2}{\mult(I)\cdot\big(\mult(I)-1\big)}
       \leq 2\cdot (k-1+\alpha)^2,
     \end{equation}
     while $\dim_\C(R/I)\leq k-1$ for $\mult(I)=1$ and the
     above inequality \eqref{eq:multpoint:1} is still satisfied. To see $\dim_\C(R/I)\leq
     k-1$ for $\mult(I)=1$ note that the ideal $I$ contains an element $g$ of order
     $1$ with $g_1=ax+by$ as homogeneous part of degree $1$ and the
     partial derivatives of $f$; applying a linear change of
     coordinates we may assume $g_1=x$ and $f=\prod_{i=1}^k
     (x-a_iy)$ with pairwise different $a_i$, and we may consider the
     negative degree lexicographical monomial ordering $>$ giving preference to
     $y$; if some $a_i=0$, then $L_>\big(\frac{\partial f}{\partial
       x}\big)=y^{k-1}$, while otherwise $L_>\big(\frac{\partial f}{\partial
       y}\big)=y^{k-1}$, so that in any case
     $\langle x,y^{k-1}\rangle\subseteq L_>(I)$, and by Proposition~
     \ref{prop:ordering} therefore
     $\dim_\C(R/I)=\dim_\C\big(R/L_>(I)\big)\leq
     \dim_\C(R/\langle x,y^{k-1}\rangle)=k-1$.
     
     Equation \eqref{eq:multpoint:1} together with Lemma
     \ref{lem:deg-fp} shows
     \begin{displaymath}
       \gamma_\alpha^{es}(M_k)\leq 2\cdot (k-1+\alpha)^2.
     \end{displaymath}
     On the other hand, considering the representative $f=x^k-y^k$, we
     have 
     \begin{displaymath}
       I^{es}(f)=\langle x^{k-1},y^{k-1}, x^ay^b\;|\;a+b=k\rangle, 
     \end{displaymath}
     and
     $I=\langle y^{k-1},x^2\rangle$ is a complete
     intersection ideal containing $I^{es}(f)$. Moreover,
     $i\big(f,x^2\big)=2k$, $\dim_\C(R/I)=2\cdot(k-1)$, thus
     \begin{displaymath}
       \gamma_\alpha^{es}(M_k)\geq\frac{\big(\alpha\cdot i(f,x^2)+(1-\alpha)\cdot\dim_\C(R/I)\big)^2}{i\big(f,x^2\big)-\dim_\C(R/I)}=
       2\cdot (k-1+\alpha)^2.
     \end{displaymath}     

     The ``in particular'' part then follows right away from Corollary
     \ref{cor:tau-ci-fat-points}. 
   \end{proof}

   Since a convenient semi-quasihomogeneous power series of
   multiplicity $2$ defines an $A_k$-singularity and one with a
   homogeneous leading term defines an ordinary multiple point, the
   following proposition together with the previous two gives upper
   bounds for all singularities defined by a convenient
   semi-quasihomogeneous representative.

   \begin{proposition}[(Semiquasihomogeneous Singularities)]\label{prop:gamma}
     Let $\ks_{p,q}$ be a singularity type with a convenient
     semi-quasihomogeneous representative $f\in R$, 
     $q>p\geq 3$. 

     Then $\gamma_\alpha^{es}(\ks_{p,q})\geq \frac{\left(q-(1-\alpha)\cdot
         \left\lfloor\frac{q}{p}\right\rfloor\right)^2}{\left\lfloor\frac{q}{p}\right\rfloor} 
     \geq \frac{q\cdot(p-1+\alpha)^2}{p}$ 
     and we obtain the following upper bound for $\gamma^{es}_\alpha(f)$:

     \begin{displaymath}
       \renewcommand{\arraystretch}{1.7}
       \begin{array}{|c||c|}
         \hline
         p,q  & \gamma_\alpha^{es}(f)\\\hline\hline
         q\geq 39 &  \leq 3\cdot(q-2+\alpha)^2 \\\hline
         \frac{q}{p}\in(1,2)  & \leq 3\cdot (q-1+\alpha)^2 \\\hline 
         \frac{q}{p}\in[2,4)  & \leq 2\cdot (q-1+\alpha)^2 \\\hline 
         \frac{q}{p}\in[4,\infty)  & \leq (q-1+\alpha)^2 \\\hline 
       \end{array}
       \renewcommand{\arraystretch}{1.4}
     \end{displaymath}
   \end{proposition}
   \begin{proof}
     To see the claimed lower bound for
     $\gamma_\alpha^{es}(\ks_{p,q})$ recall that (see
     \cite{GLS05}) 
     \begin{equation}\label{eq:eta:-1}
       I^{es}(f)=\big\langle \tfrac{\partial f}{\partial x}, 
       \tfrac{\partial f}{\partial y}, x^\alpha y^\beta\;\big|\; \alpha
       p+\beta q\geq pq\big\rangle.
     \end{equation}
     In particular, $I^{es}(f)\subseteq \big\langle
     y,x^{q-\lfloor\frac{q}{p}\rfloor}\big\rangle$,
     $\dim_\C(R/I)=q-\big\lfloor\frac{q}{p}\big\rfloor$ and
     $i(f,y)=q$, which 
     implies the claim.

     Let now $I$ be a complete intersection ideal with
     $I^{es}(f)\subseteq I$.
     Applying Lemma~\ref{lem:deg-fp} and $\degbound(I)\leq q$, we
     first of all note that 
     \begin{displaymath}
       (1+\alpha)^2\cdot \dim_\C(R/I) \leq \frac{(1+\alpha)^2\cdot(q+1)^2}{4}\leq 2\cdot(q-1+\alpha)^2.
     \end{displaymath}
     Moreover, if $\frac{q}{p}\geq 3$, then
     \begin{displaymath}
       (1+\alpha)^2\cdot \dim_\C(R/I) \leq
       \frac{(1+\alpha)^2\cdot\big(q^2+4q+3\big)}{6}\leq (q-1+\alpha)^2.
     \end{displaymath}
     since $\dim_\C(R/I)\leq \dim_\C\big(R/I^{es}(f)\big)\leq
     \frac{(p+1)\cdot(q+1)}{2}$ by \eqref{eq:eta:-1}. 

     It therefore suffices to show
     \begin{equation}
       \label{eq:eta:0}
       \lambda_\alpha(f;I,g)\leq
       \left\{
         \begin{array}{ll}
           3\cdot (q-2+\alpha)^2,& \mbox{if } q\geq 39,\\
           3\cdot (q-1+\alpha)^2,& \mbox{if } \frac{q}{p}\in (1,2),\\
           2\cdot (q-1+\alpha)^2,& \mbox{if } \frac{q}{p}\in [2,4),\\
           (q-1+\alpha)^2,& \mbox{if } \frac{q}{p}\in [4,\infty),\\
         \end{array}
       \right.
     \end{equation}
     where $g\in I$ with $i(f,g)\leq 2\cdot\dim_\C(R/I)$.
     Recall that 
     \begin{displaymath}
       \lambda_\alpha(f;I,g)=
       \frac{\big(\alpha\cdot i(f,g)+(1-\alpha)\cdot\dim_\C (R/I)\big)^2}{i(f,g)-\dim_\C(R/I)}.
     \end{displaymath}
     Fixing $I$ and considering $\lambda_\alpha(f;I,g)$ as a function in $i(f,g)$, where
     due to \eqref{eq:eta:6} the latter takes values between $\dim_\C(R/I)+1$ and
     $2\cdot\dim_\C(R/I)$, we note that the function is monotonously
     decreasing. In order to calculate an upper bound for
     $\lambda_\alpha(f;I,g)$ we may therefore replace $i(f,g)$ by some
     lower bound, which still exceeds $\dim_\C(R/I)+1$. Having done this
     we may then replace $\dim_\C(R/I)$ by an upper bound in order to
     find an upper bound for $\lambda(f;I,g)$.

     Note that for $q\geq 39$ we have 
     \begin{equation}\label{eq:eta:00}
       \frac{54}{19}\cdot (q-1+\alpha)^2\leq 3\cdot (q-2+\alpha)^2.       
     \end{equation}
          
     Fix $I$ and $g$, and let $L_{(p,q)}(g)=x^Ay^B$ be the leading term of $g$ w.\ r.\ t.\
     the weighted ordering $<_{(p,q)}$ (see Definition
     \ref{def:localordering}). By Remark \ref{rem:sqh} we know 
     \begin{equation}
       \label{eq:eta:1}
       i(f,g)\geq Ap+Bq. 
     \end{equation}
     Working with this lower bound for $i(f,g)$ we reduce the
     problem to find suitable upper bounds for $\dim_\C(R/I)$. 
     For this purpose we may assume that $L_{(p,q)}(g)$
     is minimal, and thus, in particular, $B\leq\mult(I)$. 

     If $A=0$, in view of
     Remark~\ref{rem:ci} we therefore have
     \begin{displaymath}
       B=\mult(I)\leq
       \frac{\degbound(I)+1}{2}\leq\frac{q+1}{2},
     \end{displaymath}
     and thus by Lemma \ref{lem:deg-fp} then
     \begin{equation}\label{eq:eta:2}
       \dim_\C(R/I)
       \leq
       B\cdot (q-B+1).
     \end{equation}
     Moreover, for $A=0$ Lemma \ref{lem:deg-bound-1} applies with
     $h=g$ and we get
     \begin{equation}
       \label{eq:eta:3}
       \dim_\C(R/I)
       \leq
       B\cdot q-1-\sum\limits_{i=1}^{B-1}\big\lfloor\tfrac{qi}{p}\big\rfloor   
       \leq
       B\cdot q-1-\left\lfloor\frac{q}{p}\right\rfloor\cdot \frac{B\cdot
         (B-1)}{2}.       
     \end{equation}

     Since $x^\alpha y^\beta\in I$ for $\alpha p+\beta q\geq pq$, we may assume  $Ap+Bq\leq
     pq$. But then, since $\dim_\C(R/I)\leq
     \dim_\C R\big/\big\langle \frac{\partial f}{\partial y},g,x^\alpha y^\beta\;|\;\alpha
     p+\beta q\geq pq\big\rangle$, we may
     apply  Lemma
     \ref{lem:deg-bound-2} with $h=\frac{\partial f}{\partial y}$ and $C=p-1$. This gives
     \begin{equation}\label{eq:eta:4}
       \dim_\C(R/I)
       \leq
         Ap+Bq-AB-\sum\limits_{i=1}^{A-1}\big\lfloor\tfrac{pi}{q}\big\rfloor
         -\sum\limits_{i=1}^{B-1}\big\lfloor\tfrac{qi}{p}\big\rfloor
         -\min\left\{A,\big\lceil\tfrac{q}{p}\big\rceil\right\},
     \end{equation}
     and if $B=0$ we get in addition
     \begin{equation}\label{eq:eta:5}
       \dim_\C(R/I)
       \leq 
       A\cdot(p-1).
     \end{equation}

     Finally note that by Lemma \ref{lem:schnittzahl}
     \begin{equation}\label{eq:eta:6}
       i(f,g)>\dim_\C(R/I).
     \end{equation}

     Let us now use the inequalities \eqref{eq:eta:00}-\eqref{eq:eta:6}
     to show \eqref{eq:eta:0}.
     For this we have to consider several cases for possible values of
     $A$ and $B$. 

     \begin{varthm-roman}[Case 1]
       $A=0$, $B\geq 1$.
     \end{varthm-roman}
     If $B=1$, then by \eqref{eq:eta:3} and \eqref{eq:eta:6} we have
     $\lambda_\alpha(f;I,g)\leq (q-1+\alpha)^2$.  
     
     We may thus assume that $B\geq 2$. By \eqref{eq:eta:1} and \eqref{eq:eta:2}
     \begin{displaymath}
       \lambda_\alpha(f;I,g)\leq \frac{B^2\cdot \big(q-(1-\alpha)\cdot(B-1)\big)^2}{B\cdot(B-1)}
       \leq 2\cdot (q-1+\alpha)^2. 
     \end{displaymath}
     If, moreover, $\frac{q}{p}\geq 3$, then we may apply
     \eqref{eq:eta:3} to find
     \begin{displaymath}
       \lambda_\alpha(f;I,g)\leq \frac{B^2\cdot
         \big(q-(1-\alpha)\cdot(B-1)\big)^2}{\big\lfloor\tfrac{q}{p}\big\rfloor\cdot\tfrac{B\cdot(B-1)}{2}+1}
       \leq (q-1+\alpha)^2.
     \end{displaymath}
     Taking \eqref{eq:eta:00} into account, this proves
     \eqref{eq:eta:0} in the case $A=0$ and $B\geq 1$. 

     \begin{varthm-roman}[Case 2]
       $A=1$, $B\geq 1$.
     \end{varthm-roman}
     From \eqref{eq:eta:4} we deduce
     \begin{displaymath}
       \dim_\C(R/I)\leq B\cdot(q-1)+(p-1)-\big\lfloor\tfrac{q}{p}\big\rfloor\cdot\tfrac{B\cdot(B-1)}{2}.
     \end{displaymath}
     Since $\frac{p-1+\alpha}{q-1+\alpha}\leq\frac{p}{q}$ we thus get
     \renewcommand{\arraystretch}{1.8}
     \begin{align*}
       \lambda_\alpha(f;I,g)&\leq
       \frac{\big(B+\frac{p-1+\alpha}{q-1+\alpha}\big)^2}{B
         +\big\lfloor\tfrac{q}{p}\big\rfloor\cdot\tfrac{B\cdot(B-1)}{2}+1}\cdot(q-1+\alpha)^2 \\
       & \leq 
       \left\{
         \begin{array}{lll}
           \frac{(B+\frac{1}{3})^2}{\frac{3B^2}{2}-\frac{B}{2}+1}\cdot (q-1+\alpha)^2&\leq (q-1+\alpha)^2, &
           \mbox{ if } \frac{q}{p}\geq 3,\\
           \frac{(B+\frac{1}{2})^2}{B^2+1}\cdot (q-1+\alpha)^2&\leq \frac{5}{4}\cdot(q-1+\alpha)^2, &
           \mbox{ if } \frac{q}{p}\geq 2,\\
           2\cdot\frac{(B+1)^2}{B^2+B+2}\cdot (q-1+\alpha)^2&\leq \frac{16}{7}\cdot(q-1+\alpha)^2, &
           \mbox{ if } \frac{q}{p}>1.
         \end{array}
       \right.
     \end{align*}
     Once more we are done, since $\frac{16}{7}\leq \frac{54}{19}$.
     \renewcommand{\arraystretch}{1.4}
     
     \begin{varthm-roman}[Case 3]
       $A\geq 2$, $B\geq 1$.
     \end{varthm-roman}
     Note that $\lfloor r\rfloor\geq r-1$ for any rational
     number $r$, and set $s=\frac{q}{p}$, then by \eqref{eq:eta:4}
     \begin{displaymath}
       \dim_\C(R/I) \leq
       Ap+Bq-(A-1)\cdot(B-1)-
       \frac{A\cdot(A-1)}{2s}-\frac{s\cdot B\cdot(B-1)}{2}-1-
       \min\big\{A,\lceil s\rceil\big\}.
     \end{displaymath}
     This amounts to
     \begin{multline*}
       \lambda_\alpha(f;I,g)\leq\\
       \frac{\Big(Ap+Bq-(1-\alpha)\cdot\big((A-1)\cdot(B-1)+\frac{A\cdot(A-1)}{2s}+\frac{s\cdot
           B\cdot(B-1)}{2}+1+\min\{A,\lceil s\rceil\}\big)\Big)^2}{
         (A-1)\cdot(B-1)+\frac{A\cdot(A-1)}{2s}+\frac{s\cdot B\cdot(B-1)}{2}+3}\\
       \leq
       \frac{\big(A\cdot(p-1+\alpha)+B\cdot(q-1+\alpha)\big)^2}{(A-1)\cdot(B-1)+
         \frac{A\cdot(A-1)}{2s}+\frac{s\cdot B\cdot(B-1)}{2}+3}
       \leq \varphi(A,B)\cdot (q-1+\alpha)^2,
     \end{multline*}
     where 
     \begin{displaymath}
       \varphi(A,B)=
       \frac{\big(\frac{A}{s}+B\big)^2}{(A-1)\cdot(B-1)+\frac{A\cdot(A-1)}{2s}+\frac{s\cdot B\cdot(B-1)}{2}+3}.
     \end{displaymath}
     For the last inequality we just note again that
     $\frac{p-1+\alpha}{q-1+\alpha}\leq\frac{p}{q}=\frac{1}{s}$, while for the second inequality
     a number of different cases has to be considered. We postpone
     this for a moment.

     In order to show \eqref{eq:eta:0} in the case $A\geq 2$ and
     $B\geq 1$ it now suffices to show
     \begin{equation}\label{eq:eta:8}
       \varphi(A,B)\leq
       \left\{
         \begin{array}{ll}
           \frac{54}{19},& \mbox{ if } s\geq 1,\\
           2,& \mbox{ if } s\geq 2,\\
           1,& \mbox{ if } s\geq 4.
         \end{array}
       \right.
     \end{equation}
     Elementary calculus shows that for $B\geq 1$ fixed the function
     $[2,\infty)\rightarrow \R:A\mapsto\varphi(A,B)$ takes its maximum at 
     \begin{displaymath}
       A=\max\left\{2,\frac{16-3B}{2+\frac{1}{s}}\right\}.
     \end{displaymath}

     If $B\leq 3$, then the maximum is attained at
     $A=\frac{16-3B}{2+\frac{1}{s}}$, and
     \begin{displaymath}
       \varphi(A,B)\leq
       \varphi\left(\frac{16-3B}{2+\frac{1}{s}},B\right)
       = \frac{8sB-8B+64}{4s^2B-4s^2-4sB+28s-1}.
     \end{displaymath}
     Again elementary calculus shows that the function $B\mapsto
     \varphi\left(\frac{16-3B}{2+\frac{1}{s}},B\right)$ is
     monotonously decreasing on $[1,3]$ and, therefore, 
     \begin{displaymath}
       \varphi(A,B)\leq \varphi\left(\frac{13}{2+\frac{1}{s}},1\right)
       =\frac{8s+56}{24s-1}=:\psi_1(s).
     \end{displaymath}
     Since also the function $\psi_1$ is
     monotonously decreasing on $[1,\infty)$ and
     $\psi_1(1)=\frac{64}{23}\leq \frac{54}{19}$, $\psi_1(2)=\frac{72}{47}\leq 2$
     and $\psi_1(4)=\frac{88}{95}\leq 1$ Equation \eqref{eq:eta:8} follows
     in this case.

     As soon as $B\geq 4$ the maximum for $\varphi(A,B)$ is attained
     for  $A=2$ and 
     \begin{displaymath}
       \varphi(A,B)\leq\varphi(2,B)=\frac{2\cdot(sB+2)^2}{s^3B^2-s^3B+2s^2B+4s^2+2s}.
     \end{displaymath}
     Once more elementary calculus shows that the function
     $B\mapsto\varphi(2,B)$ is monotonously decreasing on 
     $[4,\infty)$. Thus
     \begin{displaymath}
       \varphi(A,B)\leq\varphi(2,4)=\frac{4\cdot(1+2s)^2}{6s^3+6s^2+s}=:\psi_2(s).
     \end{displaymath}
     Applying elementary calculus again, we find that the function
     $\psi_2$ is monotonously decreasing on $[1,\infty)$, so that we
     are done since $\psi_2(1)=\frac{36}{13}\leq \frac{54}{19}$,
     $\psi_2(2)=\frac{50}{37}\leq 2$ and $\psi_2(4)=\frac{81}{121}\leq
     1$.

     Let us now come back to proving the missing inequality above. We
     have to show 
     \begin{displaymath}
       A+B\leq
       (A-1)\cdot(B-1)+\frac{A\cdot(A-1)}{2s}+
       \frac{s\cdot B\cdot(B-1)}{2}+1+\min\big\{A,\lceil s\rceil\big\},
     \end{displaymath}
     or equivalently 
     \begin{displaymath}
       \frac{A\cdot(A-1)}{2s}+\frac{s\cdot B\cdot(B-1)}{2}+2+
       \min\big\{A,\lceil s\rceil\big\}+AB-2A-2B\geq0.
     \end{displaymath}
     If $B\geq 2$, then $AB\geq 2A$ and $\frac{s\cdot
       B\cdot(B-1)}{2}+2+\min\big\{A,\lceil s\rceil\big\}\geq 2B$, 
     so we are done. It remains to consider the case $B=1$, and we
     have to show
     \begin{displaymath}
       A^2-A-2sA+2s\cdot\min\big\{A,\lceil s\rceil\big\}\geq 0.
     \end{displaymath}
     If $A\leq \lceil s\rceil$ or $A=2$ this is obvious. We may thus suppose
     that $A>\lceil s\rceil$ and $A\geq 3$. Since $\frac{A^2}{3}\geq
     A$ it remains to show 
     \begin{displaymath}
       \frac{2A^2}{3}-2sA+2s\cdot\lceil s\rceil\geq 0.
     \end{displaymath}
     For this 
     \begin{displaymath}
       \frac{2A^2}{3}-2sA+2s\cdot\lceil s\rceil
       \geq
       \left\{
         \begin{array}{lll}
           \frac{2A^2}{3}-2sA &
           \geq 0, & \mbox{ if } A\geq 3s,\\
           \frac{2A^2}{3}-\frac{4sA}{3} &
           \geq 0, & \mbox{ if } 2s\leq A\leq 3s,\\
           \frac{2A^2}{3}-sA &
           \geq 0, & \mbox{ if } \frac{3s}{2}\leq A\leq 2s,\\
           \frac{2A^2}{3}-\frac{2sA}{3} &
           \geq 0, & \mbox{ if } \lceil s\rceil\leq A\leq \frac{3s}{2}.
         \end{array}
       \right.
     \end{displaymath}

     \begin{varthm-roman}[Case 4]
       $A\geq 1$, $B=0$.
     \end{varthm-roman}
     Applying \eqref{eq:eta:4} and \eqref{eq:eta:5} we get
     \begin{displaymath}
       \lambda_\alpha(f;I,g)\leq
       \left\{
         \begin{array}{ll}
           \frac{A^2\cdot (p-1+\alpha)^2}{A}
           \leq 
           \left\{
             \begin{array}{l}
               \frac{A}{s^2}\cdot(q-1+\alpha)^2 \\
               A\cdot (q-2+\alpha)^2\\
             \end{array}
           \right\}
           & \mbox{ for any } A, \;\;\;\mbox{ and }\\
           \frac{A^2\cdot
             (p-1+\alpha)^2}{\sum_{i=1}^{A-1}\lfloor\frac{pi}{q}\rfloor
             +\min\{A,\lceil\frac{q}{p}\rceil\}}
           \leq \varphi_{\nu,s}(A)\cdot (q-1+\alpha)^2, & \mbox{ if } A\geq 3,
         \end{array}
       \right.
     \end{displaymath}
     where 
     \begin{displaymath}
       \varphi_{\nu,s}(A)=\frac{\frac{A^2}{s^2}}{\frac{A\cdot (A-1)}{2s}-
         (A-1)+\nu}
       =
       \frac{2A^2}{sA^2-(2s^2+s)\cdot A+2\cdot(\nu+1)\cdot s^2}
     \end{displaymath}
     with $\nu=2$ for $s\in(1,2]$ and $\nu=3$ for $s\in(2,\infty)$.

     In particular, due to the first two inequalities we may thus assume
     that 
     \begin{displaymath}
       A>
       \left\{
         \begin{array}{ll}
           3, & \mbox{ if } q\geq 39,\\
           3s^2, &\mbox{ if }s\in (1,2),\\
           2s^2, &\mbox{ if }s\in [2,4),\\
           s^2, &\mbox{ if }s\in [4,\infty).
         \end{array}
       \right.
     \end{displaymath}

     Note that $\varphi_{3,s}(A)\leq 1$ for $s\geq 4$, since 
     \begin{displaymath}
       A\geq s^2=\frac{9s^2}{16}+\frac{7s^2}{16}\geq\frac{s\cdot(1+2s)}{2\cdot(s-2)}+\frac{s}{s-2}\cdot
       \sqrt{s^2-3s+\tfrac{33}{4}}.
     \end{displaymath}
     This gives \eqref{eq:eta:0} for $s\geq 4$. 

     If now $s\in(2,4)$, then $\varphi_{3,s}$ is monotonously decreasing on
     $\big[2s^2,\infty\big)$, as is $s\mapsto \varphi_{3,s}\big(2s^2\big)$
     on $[2,4)$, and thus
     \begin{displaymath}
       \varphi_{3,s}(A)\leq \varphi_{3,s}\big(2s^2\big)=
       \frac{4s^2}{2s^3-2s^2-s+4}\leq \frac{8}{5}\leq 2,
     \end{displaymath}
     while for $s=2$ the function $\varphi_{2,2}$ is monotonously
     decreasing on $[8,\infty)$ and thus $\varphi_{2,2}(A)\leq
     \frac{16}{9}\leq 2$.
     This finishes the case $s\in[2,4)$. 

     Let's now consider the case $s\in(1,2)$ and $q\geq 39$ parallel. Applying elementary
     calculus, we find that
     $\varphi_{2,s}$ takes its maximum on $[3,\infty)$ at
     $A=\frac{12s}{1+2s}$ and is monotonously 
     decreasing on $\big[\frac{12s}{1+2s},\infty\big)$.  Moreover, the function
     $s\mapsto\varphi_{2,s}\big(\frac{12s}{1+2s}\big)$ is monotonously decreasing
     on $(1,2)$. If $s\geq \frac{7}{6}$, then 
     \begin{displaymath}
       \varphi_{2,s}(A)\leq \varphi_{2,s}\big(\tfrac{12s}{1+2s}\big)\leq
       \varphi_{2,\frac{7}{6}}\big(\tfrac{21}{5}\big)= \frac{54}{19}. 
     \end{displaymath} 
     Due to \eqref{eq:eta:00} it thus remains to consider the case
     $s\in\big(1,\frac{7}{6}\big)$ and $A>3$. If
     $A\geq 8$, then
     \begin{displaymath}
       \varphi_{2,s}(A)\leq \varphi_{2,1}(8)=\frac{64}{23}\leq\frac{54}{19},
     \end{displaymath}
     since the function
     $s\mapsto \varphi_{2,s}(8)$ is monotonously decreasing on $[1,2)$.

     So, we are finally stuck with the case $A\in\{4,5,6,7\}$ and
     $1\leq\frac{q}{p}=s\leq\frac{7}{6}$.
     We want to apply Lemma \ref{lem:deg-fp}.
     For this we note first that by Lemma \ref{lem:degb} in
     our situation $\degbound(I)\leq p+1$ and $A=\mult(I)\leq
     \frac{p+2}{2}$. But then
     \begin{displaymath}
       \dim_\C(R/I)\leq A\cdot (p-A+2)
     \end{displaymath}
     and thus, 
     \begin{displaymath}
       \lambda_\alpha(f;I,g)\leq \frac{A^2\cdot \big(p-(1-\alpha)\cdot(A-2)\big)^2}{A\cdot (A-2)}
       \leq \frac{A}{(A-2)}\cdot (q-2+\alpha)^2\leq 2\cdot (q-2+\alpha)^2.
     \end{displaymath}
     This finishes the proof.
   \end{proof}

   \begin{remark}
     In the proof of the previous proposition we achieved for almost all
     cases $\lambda_\alpha(f;I,g)\leq \frac{54}{19}\cdot(q-1+\alpha)^2$, apart from
     the single case $L_{<_{(p,q)}}(g)=x^3$. The following example
     shows that indeed in this case we cannot, in general, expect any
     better coefficient than $3$. More precisely, the example shows
     that the bound 
     \begin{displaymath}
       3\cdot (q-2+\alpha)^2
     \end{displaymath}
     is sharp for the family of singularities given by $x^q-y^{q-1}$,
     $q\geq 39$. A closer investigation should allow to lower the
     bound on $q$, but we cannot get this for all $q\geq 4$, as the
     example of $E_6$ and $E_8$ show. 

     Moreover, we give series of examples for which the bound
     $(q-1+\alpha)^2$ is sharp, respectively for which $2\cdot (q-1+\alpha)^2$ is a
     lower bound. 
   \end{remark}

   \begin{example}
     Throughout these examples $q>p\geq 3$ are integers.
     \begin{enumerate}
     \item Let $f=x^q-y^{q-1}$, then $\gamma_\alpha^{es}(f)\geq 3\cdot
       (q-2+\alpha)^2$. In particular, for $q\geq 39$,
       \begin{displaymath}
         \gamma_\alpha^{es}(f)= 3\cdot (q-2+\alpha)^2.
       \end{displaymath}
       \aus{
         For this we note that $I=\langle x^3,y^{q-2}\rangle$ is a
         complete intersection ideal in $R$ with $I^{es}(f)=\big\langle
         x^{q-1}, y^{q-2}, x^\alpha y^\beta \;\big|\; \alpha\cdot
         (q-1)+\beta q\geq q\cdot (q-1)\big\rangle\subseteq I$, since
         $2\cdot (q-1)+(q-3)\cdot q = q^2-q-2<q\cdot(q-1)$ and thus
         $x^2y^{q-3}\not\in I^{es}(f)$. This also shows that the monomial
         $x^iy^j$ with $0\leq i\leq 2$ and $0\leq j\leq q-3$ form a
         $\C$-basis of $R/I$, so that $\dim_\C(R/I)=3q-6$. Since
         $i\big(f,x^3\big)=3q-3$, the claim follows.
         }{}
     \item Let $\frac{q}{p}<2$ and $f=x^q-y^p$, then 
       \begin{displaymath}
         \gamma_\alpha^{es}(f)\geq 2\cdot (q-1+\alpha)^2.
       \end{displaymath}
       \aus{
         By the assumption on $p$ and $q$ we have $(q-2)\cdot p+q<pq$
         and hence $x^{q-2}y\not\in I^{es}(f)$. Thus $I^{es}(f)=\big\langle
         x^{q-1}, y^{p-1}, x^\alpha y^\beta \;\big|\; \alpha p
         +\beta q\geq pq\big\rangle\subseteq
         I=\langle y^2,x^{q-1}\rangle$, and we are done since
         $\dim_\C(R/I)=2q-2$ and $i\big(f,y^2\big)=2q$.
         }{}
     \item Let $f\in R$
       be convenient, semi-quasihomogeneous of $\ord_{(p,q)}(f)=pq$,
       and suppose that in $f$ no monomial $x^ky$, $k\leq
       q-2$, occurs (e.\ g.\ $f=x^q-y^p$),  then $\gamma_\alpha^{es}(f)\geq
       (q-1+\alpha)^2$. In particular, if $\frac{q}{p}\geq 4$, then
       \begin{displaymath}
         \gamma_\alpha^{es}(f)= (q-1+\alpha)^2.
       \end{displaymath}
       \aus{
         By the assumption, $I^{es}(f)\subseteq I=\langle
         x^{q-1},y\rangle$, since $\frac{\partial f}{\partial
           x}\equiv x^{q-1}\cdot u(x)\; (\mod y)$ for a unit $u$ and $\frac{\partial f}{\partial
           y}\equiv 0 \;\big(\mod \langle
         y,x^{q-1}\rangle\big)$. Hence we are done since 
         $\dim_\C(R/I)=q-1$ and $i(f,y)=q$. 
         }{}
     \item Let $f=y^3-3x^8y+3x^{12}$, then $f$ does not satisfy the
       assumptions of (c), but still
       $\gamma_\alpha^{es}(f)=(11+\alpha)^2=(q-1+\alpha)^2$. 
       \aus{
         
         For this note that $I=\langle y-x^4,x^{11}\rangle$
         contains $I^{es}(f)$, $\dim_\C(R/I)=11$ and
         $i\big(f,y-x^4\big)=12$. 
         }{}
     \item Let $f=7y^3+15x^7-21x^5y$, then $f$ is semi-quasihomogeneous
       with weights $(p,q)=(3,7)$ and convenient, but
       $\gamma_0^{es}(f)\leq 25<36=(q-1)^2$. This shows that $(q-1)^2$ is
       not a general lower bound for $\gamma_0^{es}(\ks_{p,q})$. 
       \aus{
         
         We note first that $I^{es}(f)=\langle
         x^7,y^2-x^5,x^6-x^4y\rangle$ is not a complete intersection
         and $\dim_\C\big(R/I^{es}(f)\big)=11$. Let now $I$ be a
         complete intersection ideal with $I^{es}(f)\subset I$ and let
         $h\in I$ such that $L_{<_{(3,7)}}(h)=x^Ay^B$ is minimal, in particular, 
         $\ord_{(3,7)}(h)=3A+7B$ is minimal. Then $\dim_\C(R/I)\leq 10$
         and $i(f,g)\geq 3A+7B$ for all $g\in I$.
         
         If, therefore, $3A+7B\geq 14$, then 
         \begin{displaymath}
           \frac{\dim_\C(R/I)^2}{i(f,g)-\dim_\C(R/I)}\leq 25.
         \end{displaymath}
         We may thus assume that $3A+7B\leq 13$, in particular $B<2$. If
         $B=0$, and hence $A\leq 4$, then by Lemma \ref{lem:deg-bound-2}
         $\dim_\C(R/I)\leq 2A$, so that
         \begin{displaymath}
           \frac{\dim_\C(R/I)^2}{i(f,g)-\dim_\C(R/I)}\leq 4A\leq 16.
         \end{displaymath}
         Similarly, if $B=1$ and $A=2$, then by the same Lemma
         $\dim_\C(R/I)\leq 9$ and $i(f,g)\geq 13$, so that 
         \begin{displaymath}
           \frac{\dim_\C(R/I)^2}{i(f,g)-\dim_\C(R/I)}\leq \frac{81}{4}.
         \end{displaymath}
         So it remains to consider the case $B=1$ and
         $A\in\{0,1\}$. That is $h=x^Ay+h'$ with $\ord_{(3,7)}(h')\geq
         9+3A$. Consider the ideal $J=\big\langle x^\alpha
         y^\beta\;\big|\;3\alpha+7\beta\geq 21\big\rangle\subseteq
         I$. Then $x^{4-A}\cdot h\equiv x^4y \;(\mod J)$, and thus
         $x^6-x^4y\equiv x^6\;(\mod \langle h\rangle+J)$, i.\ e.\
         $\langle h,x^6-x^4y\rangle +J=\langle h,x^6\rangle+J$. Moreover,
         $x^6\not\in \langle h\rangle+J$, so that
         $\dim_\C\big(R\big/\langle g,x^6-x^4y\rangle+J\big)\leq
         6+A$. If we can show that $\langle
         g,x^6-x^4y\rangle+J\subsetneqq I$, then
         \begin{displaymath}
           \frac{\dim_\C(R/I)^2}{i(f,g)-\dim_\C(R/I)}\leq
           \frac{(5+A)^2}{3A+7-5-A}\leq \frac{25}{2}.
         \end{displaymath}
         We are therefore done, once we know that $y^2-x^5\not\in
         \langle g,x^6\rangle +J$. Suppose there was a $g$ such
         that $gh=y^2-x^5 \;\big(\mod \langle x^6\rangle
         +J\big)$. Then $y^2=L_{<_{(3,7)}}(g)\cdot L_{<_{(3,7)}}(h)$,
         which in particular means $A=0$ and
         $L_{<_{(3,7)}}(h)=L_{<_{(3,7)}}(g)=y$. But then the coefficients
         of $1$, $x$ and $x^2$ in $h$ and $g$ must be zero, so that
         $x^5$ cannot occur with a non-zero coefficient in the
         product. This gives the desired contradiction.
         }{}
     \end{enumerate}
   \end{example}


   \section{Local Monomial Orderings}\label{sec:monomialorderings}

   Throughout the proofs of the auxilary statements in Section
   \ref{sec:sqhs} we  make use of some results from
   computer algebra concerning properties of local monomial
   orderings. In this section we recall the relevant definitions and
   results. 

   \begin{definition}\label{def:localordering}
     A \emph{monomial ordering} is a total ordering $<$ on the set
     of monomials $\left\{x^\alpha y^\beta\;\big|\;\alpha,\beta\geq
       0\right\}$ such that for all
     $\alpha,\beta,\gamma,\delta,\mu,\nu\geq 0$
     \begin{displaymath}
       x^\alpha y^\beta<x^{\gamma}y^{\delta}
       \;\;\;\Longrightarrow\;\;\;
       x^{\alpha+\mu} y^{\beta+\nu}<x^{\gamma+\mu}y^{\delta+\nu}.
     \end{displaymath}
     A monomial ordering $<$ is called \emph{local} if $1>x^\alpha
     y^\beta$ for all $(\alpha,\beta)\not=(0,0)$, and it is a local
     \emph{degree ordering} if 
     \begin{displaymath}
       \alpha+\beta>\gamma+\delta
       \;\;\;\Longrightarrow\;\;\;
       x^\alpha y^\beta<x^{\gamma}y^{\delta}.
     \end{displaymath}
     Finally, if $<$ is any local monomial ordering, then
     we define the \emph{leading monomial} $L_<(f)$ with 
     respect to $<$ of a non-zero power series $f\in R$ to be the
     maximal monomial $x^\alpha y^\beta$ such that the coefficient of $x^\alpha
     y^\beta$ in $f$ does not vanish. For $f=0$, we set $L_<(f):=0$.
     \\
     If $I\unlhd R$ is an ideal in $R$, then $L_<(I)=\langle
     L_<(f)\;|\;f\in I\rangle$ is called its \emph{leading ideal}.
   \end{definition}

   We will give now some examples of local monomial orderings which
   are used in the proofs.

   \begin{example}
     Let $\alpha,\beta,\gamma,\delta\geq 0$ be integers.
     \begin{enumerate}
     \item The \emph{negative lexicographical ordering}
       $<_{ls}$ is defined by the relation
       \begin{displaymath}
         x^\alpha y^\beta <_{ls} x^\gamma y^\delta 
         \;\;\;:\Longleftrightarrow\;\;\;
         \alpha > \gamma   \mbox{ or } (\alpha  =
         \gamma   \mbox{ and } \beta > \delta).
       \end{displaymath}
     \item The \emph{negative degree reverse lexicographical ordering}
       $<_{ds}$ is defined by the relation
       \begin{displaymath}
         x^\alpha y^\beta <_{ds} x^\gamma y^\delta 
         \;\;\;:\Longleftrightarrow\;\;\;
         \alpha +\beta  > \gamma +\delta  \mbox{ or } (\alpha +\beta  =
         \gamma +\delta  \mbox{ and } \beta > \delta).
       \end{displaymath}
     \item If positive integers $p$ and $q$ are given, then we define
       the   \emph{local weighted degree ordering}
       $<_{(p,q)}$  with weights $(p,q)$ by the relation
       \begin{displaymath}
         x^\alpha y^\beta <_{(p,q)} x^\gamma y^\delta 
         \;\;\;:\Longleftrightarrow\;\;\;
         \begin{array}[t]{l}
           \alpha p+\beta q > \gamma p+\delta q \mbox{ or } \\
           (\alpha p+\beta q =
           \gamma p+\delta q \mbox{ and } \beta < \delta).
         \end{array}
       \end{displaymath}       
     \end{enumerate}
     We note that $<_{ds}$ is a local degree ordering, while $<_{ls}$
     is not and $<_{(p,q)}$ is if and only if $p=q$.
   \end{example}

   Let us finally recall some useful properties of local orderings
   (see e.\ g.\ \cite{GP02} Corollary 7.5.6 and Proposition 5.5.7).

   \begin{proposition}\label{prop:ordering}
     Let $<$ be any local monomial ordering, and let $I$ be a
     zero-dimensional ideal in $R$.
     \begin{enumerate}
     \item The monomials of $R/L_<(I)$ form a $\C$-basis of $R/I$. In
       particular
       \begin{displaymath}
         \dim_\C(R/I)=\dim_\C\big(R/L_<(I)\big).
       \end{displaymath}
     \item If $<$ is a degree ordering, then  
       the Hilbert Samuel functions of $R/I$ and of $R/L_<(I)$
       coincide (see Definition \ref{def:hilbertsamuel}, and see also
       Remark \ref{rem:degreeordering}).
     \end{enumerate}
   \end{proposition}


   \section{The Hilbert Samuel Function}

   A useful tool in the study of the degree of zero-dimensional
   schemes and their subschemes is the Hilbert Samuel function of the
   structure sheaf, that is of the corresponding Artinian ring. 

   \begin{definition}\label{def:hilbertsamuel}
     Let $I\lhd R$ be a zero-dimensional ideal.
     \begin{enumerate}
     \item The function
       \begin{displaymath}
         H^1_{R/I}:\Z\rightarrow\Z:d\mapsto \left\{
         \begin{array}{ll}
           \dim_\C\big(R \big/ (I+\m^{d+1})\big), & d\geq 0,\\
           0, & d< 0,
         \end{array}
         \right.
       \end{displaymath}
       is called the \emph{Hilbert Samuel function} of $R/I$. 
     \item We define the \emph{slope} of the Hilbert Samuel function
       of $R/I$ to be the function
       \begin{displaymath}
         H^0_{R/I}:\N\rightarrow\N:d\mapsto
         H^1_{R/I}(d)-H^1_{R/I}(d-1).
       \end{displaymath}
       Thus
       \begin{displaymath}
         H^0_{R/I}(d)
         =\dim_\C\big(\m^d \big/ ((I\cap \m^d)+\m^{d+1})\big),
       \end{displaymath}
       is just the number $d+1$ of linearly independent
       monomials of degree $d$ in $\m^d$, minus the
       number of linearly independent monomials of degree $d$ in
       $\big(I\cap \m^d\big)+\m^{d+1}$. 
       \aus{
       
         Note that if $\overline{\m}=\m/I$ denotes the maximal ideal of
         $R/I$ and $\Gr_\m(R/I)=\bigoplus_{d\geq
           0}\overline{\m}^d/\overline{\m}^{d+1}$ the associated graded
         ring, then 
         \begin{displaymath}
           H^0_{R/I}(d)=\dim_\C\big(\overline{\m}^d/\overline{\m}^{d+1}\big)
         \end{displaymath}
         is just the dimension of the graded piece of degree $d$ of
         $\Gr_\m(R/I)$. 
         }{}
     \item Finally, we define the \emph{multiplicity} of $I$ to be 
       \begin{displaymath}
         \mult(I):=\min\big\{\mult(f)\;\big|\;0\not=f\in I\big\},
       \end{displaymath}
       and the \emph{degree bound} of $I$ as
       \begin{displaymath}
         \degbound(I):=\min\big\{d\in\N\;\big|\;\m^d\subseteq I\big\}.
       \end{displaymath}
     \end{enumerate}
   \end{definition}

   Let us gather some straight forward properties of the slope of the
   Hilbert Samuel function.

   \begin{lemma}\label{lem:hilbert-samuel}
     Let $J\subseteq I\lhd R$ be  zero-dimensional ideals.
     \begin{enumerate}
     \item $H^0_{R/I}(d)=d+1$ for all $0\leq d<\mult(I)$.
     \item $H^0_{R/I}(d)\leq H^0_{R/I}(d-1)$ for all $d\geq\mult(I)$.
     \item $H^0_{R/I}(d)\leq \mult(I)$.
     \item $H^0_{R/I}(d)=0$ for all
       $d\geq \degbound(I)$ and $H^0_{R/I}\not=0$ for all
       $d<\degbound(I)$. In particular
       \begin{displaymath}
         \dim_\C(R/I)=\sum_{d=0}^{\degbound(I)-1} H^0_{R/I}(d).
       \end{displaymath}
     \item $H^0_{R/I}(d)\leq H^0_{R/J}(d)$ for all $d\in\N$.
     \item $\degbound(I)$ and $\mult(I)$ are completely determined by
       $H^0_{R/I}$. 
     \end{enumerate}
   \end{lemma}
   \begin{proof}
     For  (a) we note that $I\subseteq\m^d$
     for all $d\leq\mult(I)$ and thus
     $H^0_{R/I}(d)=\dim_\C\big(\m^d/\m^{d+1}\big)=d+1$ for all $0\leq
     d<\mult(I)$. 

     By definition we see that $H^0_{R/I}(d)$ is just the number of
     linearly independent monomials of degree $d$ in $\m^d$, which is $d+1$, minus
     the number of linearly independent monomials, say $m_1,\ldots,m_r$, of degree $d$ in
     $\big(I\cap \m^d\big)+\m^{d+1}$. We note that then
     the set 
     \begin{displaymath}
       \{xm_1,\ldots,xm_r,ym_1,\ldots,ym_r\}\subseteq \m\cdot
       \big((I\cap \m^d)+\m^{d+1}\big) \subseteq \big(I\cap \m^{d+1}\big)+\m^{d+2}
     \end{displaymath}
     contains at least $r+1$ linearly independent monomials of degree
     $d+1$, once $r$ was non-zero. However, for
     $d=\mult(I)$ and  $g=g_{d}+h.o.t\in I$
     with homogeneous part  $g_d\not=0$ of degree $d$, we have $g_{d}\in
     \big(I\cap \m^d\big)+\m^{d+1}$, that 
     is, $d=\mult(I)$ is the smallest integer $d$ for which there is a
     monomial of degree $d$ in $\big(I\cap \m^d\big)+\m^{d+1}$.
     Thus for $d\geq \mult(I)-1$ 
     \begin{displaymath}
       H^0_{R/I}(d+1)\leq (d+2)-(r+1)=d+1-r=H^0_{R/I}(d),
     \end{displaymath}
     which proves (b), while  (c) is an immediate consequence of (a) and
     (b).

     If $d\geq \degbound(I)$, then $H^1_{R/I}(d)=\dim_\C(R/I)$ is independent
     of $d$, and hence $H^0_{R/I}(d)=0$ for all $d\geq
     \degbound(I)$. 
     In particular, 
     \begin{displaymath}
       \sum_{i=0}^{\degbound(I)-1}H^0_{R/I}(d)=H^1_{R/I}(\degbound(I)-1)-H^1_{R/I}(-1)=
       \dim_\C(R/I).
     \end{displaymath}
     Moreover, $\m^{\degbound(I)-1}+I\not=I=I+\m^{\degbound(I)}$, so
     that $H^0_{R/I}\big(\degbound(I)-1\big)\not=0$, and by (b) then
     $H^0_{R/I}(d)\not=0$ for all $d<\degbound(I)$.
     This proves  (d), and  (e) and (f) are obvious. 
   \end{proof}

   \begin{remark}\label{rem:degreeordering}
     Let $<$ be a local degree ordering on $R$, then the Hilbert
     Samuel functions of $R/I$ and of $R/L_<(I)$ coincide by
     Proposition \ref{prop:ordering}, and hence we have as well
     \begin{displaymath}
       H^0_{R/I}=H^0_{R/L_<(I)},\;\;
       \degbound(I)=\degbound\big(L_<(I)\big),\;\;
       \mbox{ and }\;\;
       \mult(I)=\mult\big(L_<(I)\big),
     \end{displaymath}
     since by the previous lemma the multiplicity and the degree bound
     only depend on the slope of the Hilbert Samuel function.
   \end{remark}

   \begin{remark}\label{rem:hilbert-samuel}
     The slope of the Hilbert Samuel function of $R/I$ gives rise to a histogram
     as the graph of the function $H^0_{R/I}$. By the Lemma
     \ref{lem:hilbert-samuel}  we
     know that up to $\mult(I)-1$ the histogram is just a staircase with
     steps of height one, and from $\mult(I)-1$ on it can only go down,
     which it eventually will do until it reaches the value
     zero for $d=\degbound(I)$.
     This means that we get a histogram of form shown in Figure \ref{fig:generalhistogram}.
     \begin{figure}[h]
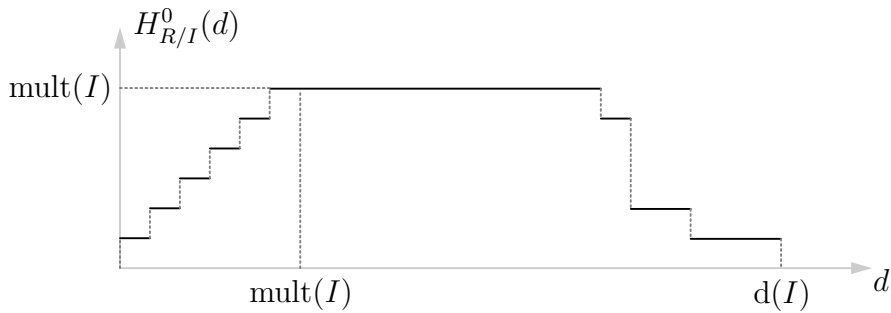

       \bigskip
       \begin{center}
         \begin{texdraw}
           \drawdim cm
           \setunitscale 0.4
           \setgray 0.8
           \arrowheadtype t:F 
           \arrowheadsize l:0.7 w:0.4
           \move (0 0) \avec (0 8) \rmove (0.4 0) \textref h:L v:C \htext{$H^0_{R/I}(d)$}
           \move (0 0) \avec (25 0) \textref h:L v:T \htext{$d$}
           \move (0 0)
           \punkte \rlvec (0 1)  \strich \rlvec (1 0)
           \punkte \rlvec (0 1)  \strich \rlvec (1 0)
           \punkte \rlvec (0 1)  \strich \rlvec (1 0)
           \punkte \rlvec (0 1)  \strich \rlvec (1 0)
           \punkte \rlvec (0 1)  \strich \rlvec (1 0)
           \punkte \rlvec (0 1)  \strich \rlvec (11 0)
           \punkte \rlvec (0 -1) \strich \rlvec (1 0)
           \punkte \rlvec (0 -3) \strich \rlvec (2 0)
           \punkte \rlvec (0 -1) \strich \rlvec (3 0)
           \punkte \rlvec (0 -1)
           \rmove (0 -0.3) \textref h:C v:T \htext{$\degbound(I)$}
           \move (6 -0.3)  \textref h:C v:T \htext{$\mult(I)$}
           \move (6 0)     \punkte \rlvec (0 5.9)
           \move (-0.3 6)  \textref h:R v:C \htext{$\mult(I)$}
           \move (0 6)     \punkte \rlvec (5 0)
         \end{texdraw}
         \medskip
         \caption{The histogram of $H^0_{R/I}$ for a general ideal $I$.}
         \label{fig:generalhistogram}
       \end{center}
     \end{figure}

     Note also, that by Lemma \ref{lem:hilbert-samuel} (a) the area of
     the histogram is just $\dim_\C(R/I)$!
   \end{remark}

   \begin{example}\label{ex:simple}
     In order to understand the slope of the Hilbert Samuel function
     better, let us consider some examples.
     \begin{enumerate}
     \item Let $f=x^2-y^{k+1}$, $k\geq 1$, and let $I=I^{ea}(f)=\langle
       x,y^k\rangle$ the equisingularity ideal of an
       $A_k$-singularity. Then $\degbound(I)=k$, $\mult(I)=1$ and
       $\dim_\C(R/I)=k$. 
       \begin{figure}[h]
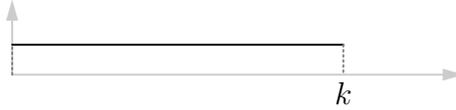

         \begin{center}
           \begin{texdraw}
             \drawdim cm
             \setunitscale 0.4
             \setgray 0.8
             \arrowheadtype t:F 
             \arrowheadsize l:0.7 w:0.4
             \move (0 0) \avec (0 2.5) 
             \move (0 0) \avec (15 0) 
             \move (0 0)
             \punkte \rlvec (0 1)  \strich \rlvec (11 0)
             \punkte \rlvec (0 -1)
             \rmove (0 -0.3) \textref h:C v:T \htext{$k$}
           \end{texdraw}
           \medskip           
           \caption{The histogram of $H^0_{R/I}$ for an $A_k$-singularity}
           \label{fig:ak-histogram}
         \end{center}
       \end{figure}
     \item Let $f=x^2y-y^{k-1}$, $k\geq 4$, and let $I=I^{ea}(f)=\langle
       xy,x^2-{\scriptsize (k-1)}\cdot y^{k-2}\rangle$ the equisingularity ideal of a
       $D_k$-singularity. Then $x^3,xy,y^{k-1}\in I$, and thus
       $\m^{k-1}\subset I$, which gives $\degbound(I)=k-1$, $\mult(I)=2$ and
       $\dim_\C(R/I)=k$, which shows that the bound in Lemma
       \ref{lem:deg-fp} need not be obtained. 
       \begin{figure}[h]
         \begin{center}
           \begin{texdraw}
             \drawdim cm
             \setunitscale 0.4
             \setgray 0.8
             \arrowheadtype t:F 
             \arrowheadsize l:0.7 w:0.4
             \move (0 0) \avec (0 2.5) 
             \move (0 0) \avec (17 0) 
             \move (0 0)
             \punkte \rlvec (0 1)  \strich \rlvec (1 0)
             \punkte \rlvec (0 1)  \strich \rlvec (1 0)
             \punkte \rlvec (0 -1) \strich \rlvec (11 0)
             \punkte \rlvec (0 -1) 
             \rmove (0 -0.3) \textref h:C v:T \htext{$k-1$}
           \end{texdraw}
           \medskip           
           \caption{The histogram of $H^0_{R/I}$ for a $D_k$-singularity}
           \label{fig:dk-histogram}
         \end{center}
       \end{figure}
     \item Let $f=x^3-y^4$ and let $I=I^{ea}(f)=\langle
       x^2,y^3\rangle$ the equisingularity ideal of an
       $E_6$-singularity. Then $\degbound(I)=4$, $\mult(I)=2$ and
       $\dim_\C(R/I)=6$.
       \\[0.1cm]
       Let $f=x^3-xy^3$ and let $I=I^{ea}(f)=\langle
       3x^2-y^3,xy^2\rangle$ the equisingularity ideal of an
       $E_7$-singularity. Then $x^3,xy^2,y^5\in I$, and thus
       $\m^5\subset I$, which gives $\degbound(I)=5$, $\mult(I)=2$ and
       $\dim_\C(R/I)=7$.
       \\[0.1cm]
       Let $f=x^3-y^5$ and let $I=I^{ea}(f)=\langle
       x^2,y^4\rangle$ the equisingularity ideal of an
       $E_8$-singularity. Then $\degbound(I)=6$, $\mult(I)=2$ and
       $\dim_\C(R/I)=8$.
       \begin{figure}[h]
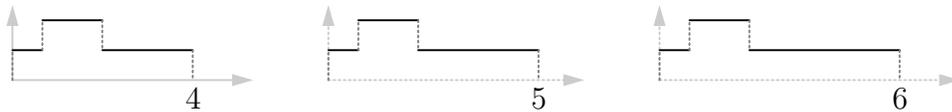

         \begin{center}
           \begin{texdraw}
             \drawdim cm
             \setunitscale 0.4
             \move (-3 0) \setgray 1 \lcir r:0.1 
             \setgray 0.8
             \arrowheadtype t:F 
             \arrowheadsize l:0.7 w:0.4
             \move (0 0) \avec (0 2.5) 
             \move (0 0) \avec (8 0) 
             \move (0 0)
             \punkte \rlvec (0 1)  \strich \rlvec (1 0)
             \punkte \rlvec (0 1)  \strich \rlvec (2 0)
             \punkte \rlvec (0 -1) \strich \rlvec (3 0)
             \punkte \rlvec (0 -1) 
             \rmove (0 -0.3) \textref h:C v:T \htext{$4$}
             \move (10.5 0)
             \setgray 0.8
             \arrowheadtype t:F 
             \arrowheadsize l:0.7 w:0.4
             \move (10.5 0) \ravec (0 2.5) 
             \move (10.5 0) \ravec (8.5 0) 
             \move (10.5 0)
             \punkte \rlvec (0 1)  \strich \rlvec (1 0)
             \punkte \rlvec (0 1)  \strich \rlvec (2 0)
             \punkte \rlvec (0 -1) \strich \rlvec (4 0)
             \punkte \rlvec (0 -1) 
             \rmove (0 -0.3) \textref h:C v:T \htext{$5$}
             \move (21.5 0)
             \setgray 0.8
             \arrowheadtype t:F 
             \arrowheadsize l:0.7 w:0.4
             \move (21.5 0) \ravec (0 2.5) 
             \move (21.5 0) \ravec (10 0) 
             \move (21.5 0)
             \punkte \rlvec (0 1)  \strich \rlvec (1 0)
             \punkte \rlvec (0 1)  \strich \rlvec (2 0)
             \punkte \rlvec (0 -1) \strich \rlvec (5 0)
             \punkte \rlvec (0 -1) 
             \rmove (0 -0.3) \textref h:C v:T \htext{$6$}
           \end{texdraw}
           \medskip           
           \caption{The histogram of $H^0_{R/I}$ for $E_6$, $E_7$ and $E_8$.}
           \label{fig:ek-histogram}
         \end{center}
       \end{figure}
     \item Let $I=\langle x^3,x^2y,y^3\rangle$, then
       $\degbound(I)=4$, $\mult(I)=3$ 
       and $\dim_\C(R/I)=7$.
       \begin{figure}[h]
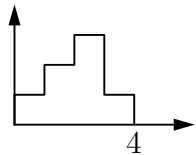

         \smallskip
         \begin{center}
           \begin{texdraw}
             \drawdim cm
             \setunitscale 0.4
             \arrowheadtype t:F 
             \arrowheadsize l:0.7 w:0.4
             \move (0 0) \avec (0 4) 
             \move (0 0) \avec (6 0) 
             \move (0 0)
             \rlvec (0 1) \rlvec (1 0)
             \rlvec (0 1) \rlvec (1 0)
             \rlvec (0 1) \rlvec (1 0)
             \rlvec (0 -2) \rlvec (1 0)
             \rlvec (0 -1) 
             \rmove (0 -0.3) \textref h:C v:T \htext{$4$}
           \end{texdraw}
           \medskip
           \caption{The histogram of $H^0_{R/I}$ for $I=\langle x^3,x^2y,y^3\rangle$.}
         \end{center}
       \end{figure}
     \end{enumerate}
   \end{example}

   The following result providing a lower bound for the minimal number
   of generators of a zero-dimensional ideal in $R$ is due to A.\
   Iarrobino. 

   \begin{lemma}\label{lem:ci}
     Let $I\lhd R$ be a zero-dimensional ideal. Then $I$
     cannot be generated by less than $1+\sup\left\{H^0_{R/I}(d-1)-H^0_{R/I}(d)\;\big|\;d\geq
       \mult(I)\right\}$ elements.

     In particular, if $I$ is a complete intersection ideal then for
     $d\geq\mult(I)$
     \begin{displaymath}
       H^0_{R/I}(d-1)-1\leq H^0_{R/I}(d)\leq H^0_{R/I}(d-1).
     \end{displaymath}
   \end{lemma}
   \begin{proof}
     See \cite{Iar77} Theorem 4.3 or \cite{Bri77} Proposition III.2.1.
   \end{proof}

   Moreover, by the Lemma of Nakayama and Proposition
   \ref{prop:ordering} we can compute the minimal number of
   generators for a zero-dimensional ideal exactly.

   \begin{lemma}\label{lem:nak}
     Let $I\lhd R$ be zero-dimensional ideal and let $<$ denote any
     local ordering on $R$. Then 
     the minimal number of generators of $I$ is
     \begin{displaymath}
       \dim_\C(I/\m I)=\dim_\C\big(R/L_<(I)\big)-\dim_\C\big(R/L_<(\m I)\big).
     \end{displaymath}
   \end{lemma}

   \begin{remark}\label{rem:ci}
     If we apply Lemma \ref{lem:ci} to a zero-dimensional complete
     intersection ideal $I\lhd R$, i.\ e.\ a zero-dimensional ideal generated
     by two elements, then we know that the histogram of $H^0_{R/I}$
     will be as shown in Figure \ref{fig:ci-histogram};
     \begin{figure}[h]
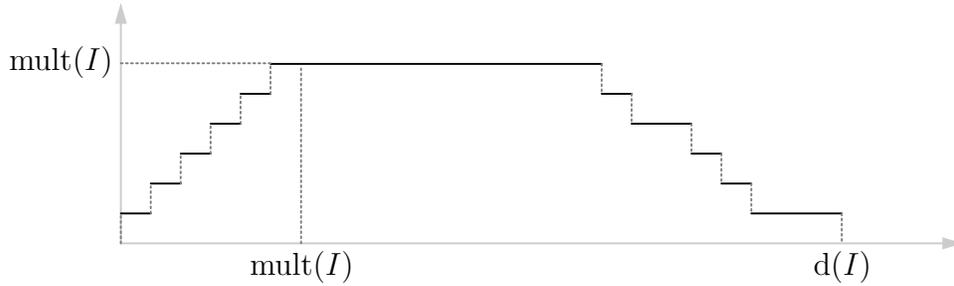

       \bigskip
       \begin{center}
         \begin{texdraw}
           \drawdim cm
           \setunitscale 0.4
           \setgray 0.8
           \arrowheadtype t:F
           \arrowheadsize l:0.7 w:0.4
           \move (0 0) \avec (0 8) \rmove (0.4 0) 
           \move (0 0) \avec (28 0) \textref h:L v:T
           \move (0 0)
           \punkte \rlvec (0 1)  \strich \rlvec (1 0)
           \punkte \rlvec (0 1)  \strich \rlvec (1 0)
           \punkte \rlvec (0 1)  \strich \rlvec (1 0)
           \punkte \rlvec (0 1)  \strich \rlvec (1 0)
           \punkte \rlvec (0 1)  \strich \rlvec (1 0)
           \punkte \rlvec (0 1)  \strich \rlvec (11 0)
           \punkte \rlvec (0 -1) \strich \rlvec (1 0)
           \punkte \rlvec (0 -1) \strich \rlvec (2 0)
           \punkte \rlvec (0 -1) \strich \rlvec (1 0)
           \punkte \rlvec (0 -1) \strich \rlvec (1 0)
           \punkte \rlvec (0 -1) \strich \rlvec (3 0)
           \punkte \rlvec (0 -1)
           \rmove (0 -0.3) \textref h:C v:T \htext{$\degbound(I)$}
           \move (6 -0.3)  \textref h:C v:T \htext{$\mult(I)$}
           \move (6 0)     \punkte \rlvec (0 5.9)
           \move (-0.3 6)  \textref h:R v:C \htext{$\mult(I)$}
           \move (0 6)     \punkte \rlvec (5 0)
         \end{texdraw}
         \medskip
         \caption{The histogram of $H^0_{R/I}$ for a complete
           intersection.}\label{fig:ci-histogram}
       \end{center}
     \end{figure}
     that is, up to the value $d=\mult(I)$ the histogram of $H^0_{R/I}$ is an
     ascending staircase with steps of height and length one, then it remains
     constant for a while, and finally it is a descending staircase
     again with steps of height one, but a possibly longer length.
     In particular we see that 
     \begin{equation}\label{eq:mult-deg}
       \mult(I)\leq
       \left\{
         \begin{array}{ll}
           \frac{\degbound(I)+1}{2}, & \mbox { if } \degbound(I)
           \mbox{ is odd},\\
           \frac{\degbound(I)}{2}, & \mbox { if } \degbound(I)
           \mbox{ is even}.
         \end{array}
       \right.
     \end{equation}
   \end{remark}

   \begin{example}
     Let $I=\m^k$ for $k\geq 1$. Then $\degbound(I)=\mult(I)=k$
     and $\dim_\C(R/I)=\binom{k+1}{2}$.
     \begin{figure}[h]
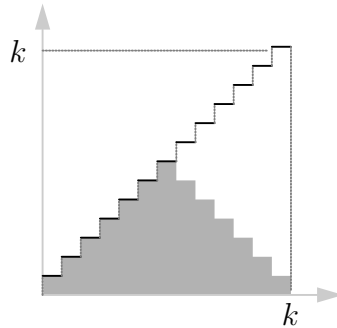

       \smallskip
       \begin{center}
         \begin{texdraw}
           \drawdim cm
           \setunitscale 0.25
           \arrowheadtype t:F 
           \arrowheadsize l:1.4 w:0.8
           \move (0 0)
           \rlvec (0 1) \rlvec (1 0)
           \rlvec (0 1) \rlvec (1 0)
           \rlvec (0 1) \rlvec (1 0)
           \rlvec (0 1) \rlvec (1 0)
           \rlvec (0 1) \rlvec (1 0)
           \rlvec (0 1) \rlvec (1 0)
           \rlvec (0 1) \rlvec (1 0)
           \rlvec (0 -1) \rlvec (1 0)
           \rlvec (0 -1) \rlvec (1 0)
           \rlvec (0 -1) \rlvec (1 0)
           \rlvec (0 -1) \rlvec (1 0)
           \rlvec (0 -1) \rlvec (1 0)
           \rlvec (0 -1) \rlvec (1 0)
           \rlvec (0 -1) \rlvec (1 0)
           \ifill f:0.7
           \setgray 0.8
           \move (0 0) \avec (0 15.5) \rmove (0.4 0) 
           \move (0 0) \avec (16 0) 
           \move (0 0)
           \punkte \rlvec (0 1)  \strich \rlvec (1 0)
           \punkte \rlvec (0 1)  \strich \rlvec (1 0)
           \punkte \rlvec (0 1)  \strich \rlvec (1 0)
           \punkte \rlvec (0 1)  \strich \rlvec (1 0)
           \punkte \rlvec (0 1)  \strich \rlvec (1 0)
           \punkte \rlvec (0 1)  \strich \rlvec (1 0)
           \punkte \rlvec (0 1)  \strich \rlvec (1 0)
           \punkte \rlvec (0 1)  \strich \rlvec (1 0)
           \punkte \rlvec (0 1)  \strich \rlvec (1 0)
           \punkte \rlvec (0 1)  \strich \rlvec (1 0)
           \punkte \rlvec (0 1)  \strich \rlvec (1 0)
           \punkte \rlvec (0 1)  \strich \rlvec (1 0)
           \punkte \rlvec (0 1)  \strich \rlvec (1 0)
           \punkte \rlvec (0 -13)
           \rmove (0 -0.6) \textref h:C v:T \htext{$k$}
           \move (-0.8 13)  \textref h:R v:C \htext{$k$}
           \move (0 13)     \punkte \rlvec (12 0)
         \end{texdraw}
         \medskip
         \caption{The histogram of $H^0_{R/\m^k}$. The shaded region
           is the maximal possible value of $\dim_\C(R/I)$ for a complete intersection
           ideal $I$ containing  $\m^k$.}
         \label{fig:fp-histogram}
       \end{center}
     \end{figure}     
   \end{example}

   \begin{lemma}\label{lem:deg-fp}
     Let $I\lhd R$ be a zero-dimensional 
     complete intersection ideal, then  
     \begin{displaymath}
       \dim_\C(R/I)\leq\big(\degbound(I)-\mult(I)+1\big)\cdot \mult(I).
     \end{displaymath}
     In particular
     \begin{displaymath}
       \dim_\C(R/I)\leq 
       \left\{
         \begin{array}{ll}
           \frac{(\degbound(I)+1)^2}{4}, & \mbox{ if } \degbound(I) \mbox{ odd},\\
           \frac{\degbound(I)^2+2\degbound(I)}{4}, & \mbox{ if } \degbound(I) \mbox{ even}.\\
         \end{array}
       \right.
     \end{displaymath}
   \end{lemma}
   \begin{proof}
     By Remark \ref{rem:hilbert-samuel} we have to find an upper bound
     for the area $A$ of the histogram of $H^0_{R/I}$. This area would
     be maximal, if in the descending part the steps had all length
     one, i.\ e.\ if the histogram was as shown in Figure \ref{fig:ci-opt-histogram}. 
     \begin{figure}[h]
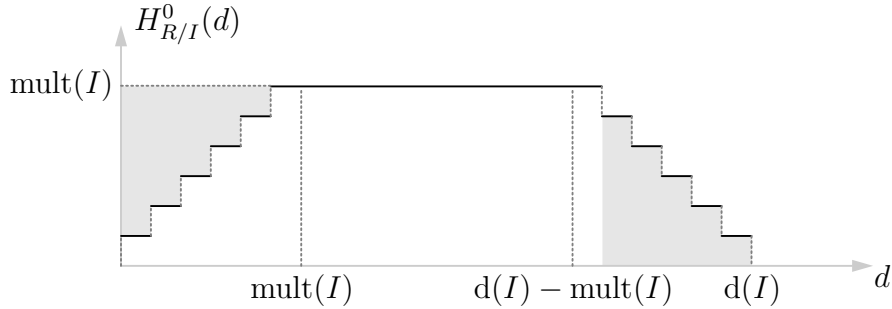

       \bigskip
       \begin{center}
         \begin{texdraw}
           \drawdim cm
           \setunitscale 0.4
           \move (16 0)
           \rlvec (0 5) \rlvec (1 0) \rlvec (0 -1) 
           \rlvec (1 0) \rlvec (0 -1) \rlvec (1 0) \rlvec (0 -1) 
           \rlvec (1 0) \rlvec (0 -1) \rlvec (1 0) \rlvec (0 -1) 
           \ifill f:0.9
           \move (0 1)  
           \rlvec (1 0) \rlvec (0 1) \rlvec (1 0) \rlvec (0 1) 
           \rlvec (1 0) \rlvec (0 1) \rlvec (1 0) \rlvec (0 1) 
           \rlvec (1 0) \rlvec (0 1) \rlvec (-5 0) \rlvec (0 -5)
           \ifill f:0.9
           \setgray 0.8
           \arrowheadtype t:F
           \arrowheadsize l:0.7 w:0.4
           \move (0 0) \avec (0 8) \rmove (0.4 0) \textref h:L v:C \htext{$H^0_{R/I}(d)$}
           \move (0 0) \avec (25 0) \textref h:L v:T \htext{$d$}
           \move (0 0)
           \punkte \rlvec (0 1)  \strich \rlvec (1 0)
           \punkte \rlvec (0 1)  \strich \rlvec (1 0)
           \punkte \rlvec (0 1)  \strich \rlvec (1 0)
           \punkte \rlvec (0 1)  \strich \rlvec (1 0)
           \punkte \rlvec (0 1)  \strich \rlvec (1 0)
           \punkte \rlvec (0 1)  \strich \rlvec (11 0)
           \punkte \rlvec (0 -1) \strich \rlvec (1 0)
           \punkte \rlvec (0 -1) \strich \rlvec (1 0)
           \punkte \rlvec (0 -1) \strich \rlvec (1 0)
           \punkte \rlvec (0 -1) \strich \rlvec (1 0)
           \punkte \rlvec (0 -1) \strich \rlvec (1 0)
           \punkte \rlvec (0 -1)
           \rmove (0 -0.3) \textref h:C v:T \htext{$\degbound(I)$}
           \move (6 -0.3)  \textref h:C v:T \htext{$\mult(I)$}
           \move (6 0)     \punkte \rlvec (0 5.9)
           \move (15 -0.3) \textref h:C v:T \htext{$\degbound(I)-\mult(I)$}
           \move (15 0)    \punkte \rlvec (0 5.9)
           \move (-0.3 6)  \textref h:R v:C \htext{$\mult(I)$}
           \move (0 6)     \punkte \rlvec (5 0)
         \end{texdraw}
         \medskip
         \caption{Maximal possible area.}
         \label{fig:ci-opt-histogram}
       \end{center}
     \end{figure}
     Since the two shaded regions  have
     the same area, we get  
     \begin{displaymath}
       A\leq\big(\degbound(I)-\mult(I)+1\big)\cdot \mult(I). 
     \end{displaymath}
     Consider now the function
     \begin{displaymath}
       \varphi:\Big[\mult(I),\tfrac{\degbound(I)+1}{2}\Big]\longrightarrow\R: x\mapsto 
       \big(\degbound(I)-x+1\big)\cdot x,  
     \end{displaymath}
     then this function is monotonously increasing, which finishes the
     proof in view of Equation \eqref{eq:mult-deg}.
   \end{proof}

   \begin{corollary}\label{cor:tau-ci-fat-points}
     For an ordinary $m$-fold point $M_m$ we have
     \begin{displaymath}
       \tau_{ci}^{es}(M_m)=
       \left\{
         \begin{array}{cl}
           \frac{(m+1)^2}{4}, & \mbox{ if } m\geq 3 \mbox{ odd},\\
           \frac{m^2+2m}{4}, & \mbox{ if } m\geq 4 \mbox{ even},\\
           1,& \mbox{ if } m=2.
         \end{array}
       \right.
     \end{displaymath}
   \end{corollary}
   \begin{proof}
     Let $f$ be a representative of $M_m$. Then 
     \begin{displaymath}
       I^{es}(f)=\left\langle
       \frac{\partial f}{\partial x},\frac{\partial f}{\partial x}\right\rangle + \m^m,
     \end{displaymath}
     and as in the proof of Proposition \ref{prop:gamma-ordinary} we may
     assume that $f$ is a homogeneous of degree $m$. 

     In particular, if $m=2$, then
     $I^{es}(f)=\m$ is a complete intersection and
     $\tau^{es}_{ci}(M_2)=1$. We may therefore assume that $m\geq 3$.

     For any complete
     intersection ideal $I$ with $\m^m\subset I^{es}(f)\subseteq I$ we
     automatically have $\degbound(I)\leq m$, and by Lemma
     \ref{lem:deg-fp} 
     \begin{displaymath}
       \tau_{ci}^{es}(f)\leq 
       \left\{
         \begin{array}{ll}
           \frac{(m+1)^2}{4}, & \mbox{ if } m \mbox{ odd},\\
           \frac{m^2+2m}{4}, & \mbox{ if } m\geq 4 \mbox{ even}.
         \end{array}
       \right.
     \end{displaymath}
     Consider now the representative $f=x^m-y^m$.
     If $m=2k$ is even, then the ideal $I= \langle
     x^k,y^{k+1}\rangle$ is a complete intersection with
     $I^{es}(f)\subset I$ and 
     \begin{displaymath}
       \tau_{ci}^{es}(f)\geq\dim_\C(R/I)=k^2+k=\frac{m^2+2m}{4}.
     \end{displaymath}
     Similarly, if $m=2k-1$ is odd, then the ideal $I= \langle
     x^k,y^k\rangle$ is a complete intersection with
     $I^{es}(f)\subset I$ and 
     \begin{displaymath}
       \tau_{ci}^{es}(f)\geq\dim_\C(R/I)=k^2=\frac{m^2+2m+1}{4}.
     \end{displaymath}     
   \end{proof}
   
   \aus{
   \begin{remark}\label{rem:non-ci}
     Let $I\lhd R$ be any zero-dimensional 
     ideal, not necessarily a complete intersection, then still  
     \begin{displaymath}
       \dim_\C(R/I)\leq\left(\degbound(I)-\frac{\mult(I)-1}{2}\right)\cdot \mult(I).
     \end{displaymath}
   \end{remark}
   \begin{proof}
     The proof is the same as for the complete intersection
     ideal, just that we cannot ensure that the histogram goes down to
     zero at $\degbound(I)$ with steps of size one. The dimension is
     thus bounded by the region of the histogram in Figure
     \ref{fig:histogram-non-ci}. 
     \begin{figure}[h]
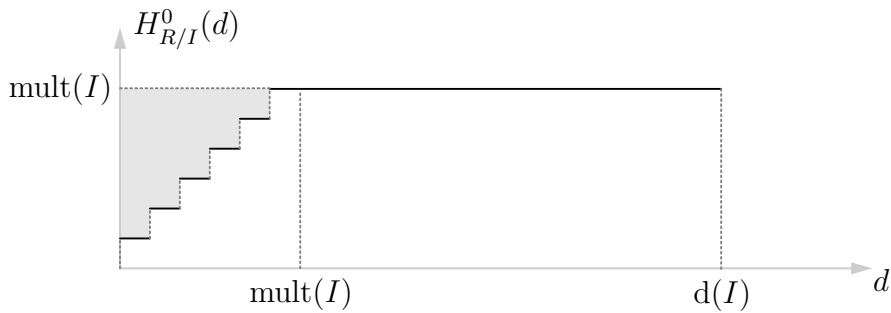

       \bigskip
       \begin{texdraw}
         \drawdim cm
         \setunitscale 0.4
         \move (0 1)  
         \rlvec (1 0) \rlvec (0 1) \rlvec (1 0) \rlvec (0 1) 
         \rlvec (1 0) \rlvec (0 1) \rlvec (1 0) \rlvec (0 1) 
         \rlvec (1 0) \rlvec (0 1) \rlvec (-5 0) \rlvec (0 -5)
         \ifill f:0.9
         \setgray 0.8
         \arrowheadtype t:F
         \arrowheadsize l:0.7 w:0.4
         \move (0 0) \avec (0 8) \rmove (0.4 0) \textref h:L v:C \htext{$H^0_{R/I}(d)$}
         \move (0 0) \avec (25 0) \textref h:L v:T \htext{$d$}
         \move (0 0)
         \punkte \rlvec (0 1)  \strich \rlvec (1 0)
         \punkte \rlvec (0 1)  \strich \rlvec (1 0)
         \punkte \rlvec (0 1)  \strich \rlvec (1 0)
         \punkte \rlvec (0 1)  \strich \rlvec (1 0)
         \punkte \rlvec (0 1)  \strich \rlvec (1 0)
         \punkte \rlvec (0 1)  \strich \rlvec (15 0)
         \punkte \rlvec (0 -6) 
         \rmove (0 -0.3) \textref h:C v:T \htext{$\degbound(I)$}
         \move (6 -0.3)  \textref h:C v:T \htext{$\mult(I)$}
         \move (6 0)     \punkte \rlvec (0 5.9)
         \move (-0.3 6)  \textref h:R v:C \htext{$\mult(I)$}
         \move (0 6)     \punkte \rlvec (5 0)
       \end{texdraw}
       \medskip
       \caption{Maximal possible area.}
       \label{fig:histogram-non-ci}
     \end{figure}
     \end{proof}
     }{}


   \section{Semi-Quasihomogeneous Singularities}\label{sec:sqhs}

   \begin{definition}\label{def:sqh}
     A non-zero polynomial of the form $f=\sum_{\alpha\cdot p+\beta\cdot
       q=d}a_{\alpha,\beta}x^\alpha y^\beta$ is called
     \emph{quasihomogeneous} of $(p,q)$-degree $d$. 
     Thus the Newton polygon of a quasihomogeneous polynomial
     has just one side of slope $-\frac{p}{q}$.

     A quasihomogeneous polynomial is said to be
     \emph{non-degenerate} if it is reduced, that is if it has no
     multiple factors, and it is said to be \emph{convenient} if
     $\frac{d}{p},\frac{d}{q}\in\Z$ and $a_{\frac{d}{p},0}$ and $a_{0,\frac{d}{q}}$ are non-zero, that is
     if the Newton polygon meets the $x$-axis and the $y$-axis. 

     If $f=f_0+f_1$ with $f_0$ quasihomogeneous of $(p,q)$-degree $d$
     and for any monomial $x^\alpha y^\beta$ occurring in $f_1$ with a
     non-zero coefficient we have $\alpha \cdot p+\beta\cdot
     q> d$, we say that $f$ is of $(p,q)$-\emph{order} $d$, and we call $f_0$
     the $(p,q)$-\emph{leading form} of $f$ and denote it by
     $\lead_{(p,q)}(f)$. We denote the $(p,q)$-order of 
     $f$ by $\ord_{(p,q)}(f)$.

     A power series $f\in R$ is said to be
     \emph{semi-quasihomogeneous} with respect to the weights $(p,q)$  if the
     $(p,q)$-leading form is non-degenerate. 
   \end{definition}

   \begin{remark}\label{rem:sqh}
     Let $f\in R$ with $\deg_{(p,q)}(f)=pq$ and let $f_0$ denote its
     $(p,q)$-leading form.
     \begin{enumerate}
     \item If $\gcd(p,q)=r$, then $f_0$ has $r$
       factors of the form $a_ix^\frac{q}{r}-b_iy^\frac{p}{r}$,
       $i=1,\ldots,r$.

       If, moreover, $f_0$ is 
       non-degenerate, then these will all be irreducible and pairwise
       different, i.\ e.\ not scalar multiples of each other.
     \item If $f$ is irreducible, then $f_0$ has only one irreducible
       factor, possibly of higher multiplicity.
     \item If $f_0$ is non-degenerate, then $f$ has $r=\gcd(p,q)$
       branches $f_1,\ldots,f_r$, which are all semi-quasihomogeneous
       with irreducible
       $(p,q)$-leading form $a_ix^\frac{q}{r}-b_iy^\frac{p}{r}$ for
       pairwise distinct points $(a_i:b_i)\in \PC^1$,
       $i=1,\ldots,r$. 

       The characteristic exponents of $f_i$ are
       $\frac{q}{r}$ and $\frac{p}{r}$ for all $i=1,\ldots,r$, and
       thus $f_i$ admits a parametrisation of the form
       \begin{displaymath}
         \big(x_i(t),y_i(t)\big)=\Big(\alpha_i
         t^\frac{p}{r}+h.o.t,\beta_i t^\frac{q}{r}+h.o.t\Big).
       \end{displaymath}
     \item If $f_0$ is non-degenerate, i.\ e.\ $f$ is semi-quasihomogeneous,
       and $g\in R$, then
       \begin{displaymath}
         i(f,g)\geq\ord_{(p,q)}(g).
       \end{displaymath}
     \end{enumerate}
   \end{remark}
   \begin{proof}\leererpunkt
     \begin{enumerate}
     \item If $\alpha p+\beta q=pq$, then $p\;|\;\beta q$ and hence
       $p\;|\;\beta r$, so that $\beta\cdot\frac{r}{p}$ is a natural
       number. Similarly $\alpha\cdot\frac{r}{q}$ is a natural
       number. We may therefore consider the transformation
       \begin{displaymath}
         f_0\big(x^\frac{r}{q},y^\frac{r}{p}\big)\in \C[x,y]_r
       \end{displaymath}
       which is a homogeneous polynomial of degree $r$. Thus
       $f_0\big(x^\frac{r}{q},y^\frac{r}{p}\big)$ factors in $r$
       linear factors $a_ix-b_iy$, $i=1,\ldots,r$, so that $f_0$
       factors as
       \begin{equation}\label{eq:newton}
         f_0=\prod_{i=1}^r\big(a_ix^\frac{q}{r}-b_iy^\frac{p}{r}\big).
       \end{equation}
       Since $\gcd\big(\frac{p}{r},\frac{q}{r}\big)=1$, the factors
       $a_ix^\frac{q}{r}-b_iy^\frac{p}{r}$ are irreducible once
       neither $a_i$ nor $b_i$ is zero. 

       If $f_0$ is non-degenerate, then the irreducible factors of
       $f_0$ are pairwise distinct. So, $a_i=0$ implies $r=p$ and
       still $a_ix^\frac{q}{r}-b_iy^\frac{p}{r}=b_iy$ irreducible,
       while $b_i=0$ similarly gives $r=q$ and
       $a_ix^\frac{q}{r}-b_iy^\frac{p}{r}=a_ix$ irreducible. Thus, in
       any case the factors in \eqref{eq:newton} are irreducible and,
       hence, pairwise distinct.
     \item With the notation from Lemma \ref{lem:brieskorn} and the
       factorisation of $f_0$ from \eqref{eq:newton} we get  
       \begin{displaymath}
         g=\frac{\prod_{i=1}^ra_i u^\frac{bq}{r}v^\frac{pq}{r^2}-b_i u^\frac{ap}{r}v^\frac{pq}{r^2}}{u^{ap}v^\frac{pq}{r}}
         =\prod_{i=1}^r(a_iu-b_i).
       \end{displaymath}
       By assumption $f$ is irreducible, hence according to Lemma
       \ref{lem:brieskorn} $g$ has at most one, possibly repeated,
       zero. But thus the factors of $f_0$ all coincide -- up to
       scalar multiple. 
     \item The first assertion is an immediate consequence from (a)
       and (b), while the ``in particular'' part follows by Puiseux
       expansion. 
     \item Let $g_0$ be the $(p,q)$-leading form of $g$. Using the
       notation from (c) we have 
       \begin{multline*}
         \;\;\;\;\;\;\;\;\;
         i(f,g)=\sum_{i=1}^r i(f_i,g)=\sum_{i=1}^r
         \ord\big(g(x_i(t),y_i(t))\big)\\
         =\sum_{i=1}^r
         \ord\Big(g_0\big(\alpha_i t^\frac{p}{r},\beta_i t^\frac{q}{r}\big)+h.o.t\Big)
         \geq\sum_{i=1}^r\frac{\ord_{(p,q)}(g)}{r}=\ord_{(p,q)}(g).
       \end{multline*}
     \end{enumerate}
   \end{proof}

   \begin{lemma}\label{lem:brieskorn}
     Let $f\in R$ with $\ord_{(p,q)}(f)=pq$ and let $f_0$ denote its
     $(p,q)$-leading form.
     Let $r=\gcd(p,q)$ and $a,b\geq 0$ such that $qb-pa=r$. Finally set 
     \begin{displaymath}
       g=\frac{f_0\big(u^bv^\frac{p}{r},u^av^\frac{q}{r}\big)}{u^{ap}v^\frac{pq}{r}}
       \in \C[u].
     \end{displaymath}
     Then the number of different zeros of $g$ is a lower bound for
     the number of branches of $f$. 
   \end{lemma}
   \begin{proof}
     See \cite{BK86} Remark on p.\ 480. 
   \end{proof}

   The following investigations are crucial for the proof of
   Proposition \ref{prop:gamma}. 

   \begin{lemma}\label{lem:deg-bound-1}
     Let $f\in R$ be convenient semi-quasihomogeneous with leading
     form $f_0$ and $\ord_{(p,q)}(f)=pq$, let $I=\big\langle x^\alpha
     y^\beta\;\big|\;\alpha p+\beta q\geq pq\big\rangle$, and let
     $h\in R$. Then  
     \begin{displaymath}
       \dim_\C R/\big(\langle h\rangle+I^{es}(f)\big)<\dim_\C R/\big(\langle h\rangle
       +I\big).
     \end{displaymath}
     In particular, if $L_{(p,q)}(h)=y^B$ with $B\leq p$, then
     \begin{displaymath}
       \dim_\C R/\langle h\rangle+I^{es}(f)\leq Bq-1-\sum_{i=1}^{B-1}\big\lfloor\tfrac{qi}{p}\big\rfloor.
     \end{displaymath}
   \end{lemma}
   \begin{proof}
     As
     \begin{displaymath}
       I^{es}(f)=\big\langle \tfrac{\partial f}{\partial x},
       \tfrac{\partial f}{\partial y}\big\rangle + I,
     \end{displaymath}
     it suffices to show that 
     \begin{displaymath}
       I^{es}(f)\not\subseteq\langle h\rangle +I,
     \end{displaymath}
     which is the same as showing that not both 
     $\frac{\partial f}{\partial x}$ and
     $\frac{\partial f}{\partial y}$ belong to $\langle h\rangle +I$. 

     Suppose the contrary, that is, there are $h_x, h_y\in R$ such
     that 
     \begin{displaymath}
       \tfrac{\partial f}{\partial x}\equiv h_x\cdot h\; (\mod I)
       \;\;\;\;\mbox{ and }\;\;\;\;
       \tfrac{\partial f}{\partial y}\equiv h_y\cdot h\; (\mod I).
     \end{displaymath}
     We note that 
     \begin{displaymath}
       \lead_{(p,q)}\big(\tfrac{\partial f}{\partial x}\big)
       =\tfrac{\partial f_0}{\partial x}
       \;\;\;\;\mbox{ and }\;\;\;\;
       \lead_{(p,q)}\big(\tfrac{\partial f}{\partial y}\big)
       =\tfrac{\partial f_0}{\partial y},
     \end{displaymath}
     and none of the monomials involved is contained in $I$. Therefore 
     \begin{displaymath}
       \lead_{(p,q)}(h_x)\cdot\lead_{(p,q)}(h)
        =\tfrac{\partial f_0}{\partial x}
       \;\;\;\;\mbox{ and }\;\;\;\;
       \lead_{(p,q)}(h_y)\cdot\lead_{(p,q)}(h)
        =\tfrac{\partial f_0}{\partial y},
     \end{displaymath}
     which in particular implies that $\tfrac{\partial f_0}{\partial
       x}$ and $\tfrac{\partial f_0}{\partial y}$ have a common
     factor. This, however,  is then a multiple factor of the quasihomogeneous polynomial
     $f_0$, in contradiction to $f$ being semi-quasihomogeneous.

     \begin{figure}[h]
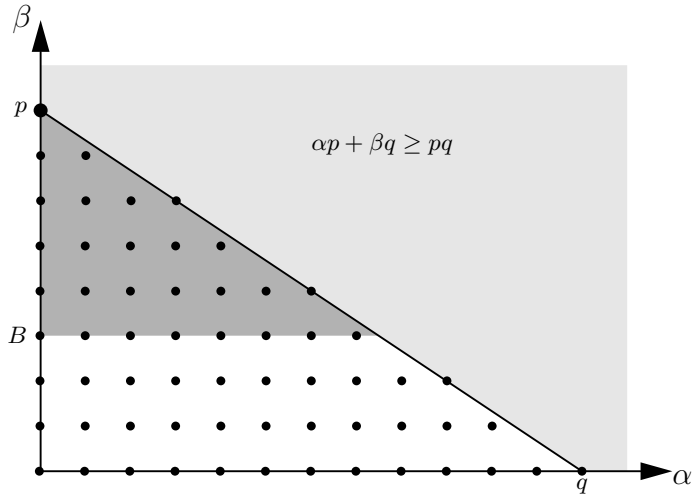

       \bigskip
       \begin{center}
         \begin{texdraw}
           \drawdim cm
           \setunitscale 0.6
           \setgray 0.8
           \arrowheadtype t:F 
           \arrowheadsize l:0.7 w:0.4
           \move (0 8) \lvec (12 0) \lvec (13 0) \lvec (13 9) \lvec (0 9) \ifill f:0.9
           \move (0 8) \lvec (0 3) \lvec (7.5 3) \ifill f:0.7
           \setgray 0
           \move (0 0) \avec (0 10) \rmove (-0.2 0)\textref h:R v:C \htext{$\beta$}
           \move (0 0) \avec (14 0) \textref h:L v:T \htext{$\alpha$}
           \move (0 8) \lvec (12 0)
           \move (0 8) \fcir f:0 r:0.15 
           \move (-0.3 8)\textref h:R v:C {\scriptsize\htext{$p$}}
           \move (-0.3 3)\textref h:R v:C {\scriptsize\htext{$B$}}
           \move (12 -0.2)\textref h:C v:T {\scriptsize\htext{$q$}}
           \move (6 7)\textref h:L v:B {\scriptsize\htext{$\alpha
               p+\beta q \geq pq$}}
           \move (1 7) \fcir f:0 r:0.1 
           \rmove (-1 0) \fcir f:0 r:0.1 
           \move (3 6) \fcir f:0 r:0.1 
           \rmove (-1 0) \fcir f:0 r:0.1 \rmove (-1 0) \fcir f:0 r:0.1
           \rmove (-1 0) \fcir f:0 r:0.1 
           \move (4 5) \fcir f:0 r:0.1 
           \rmove (-1 0) \fcir f:0 r:0.1 \rmove (-1 0) \fcir f:0 r:0.1
           \rmove (-1 0) \fcir f:0 r:0.1 \rmove (-1 0) \fcir f:0 r:0.1
           \move (6 4) \fcir f:0 r:0.1 
           \rmove (-1 0) \fcir f:0 r:0.1 \rmove (-1 0) \fcir f:0 r:0.1 
           \rmove (-1 0) \fcir f:0 r:0.1 \rmove (-1 0) \fcir f:0 r:0.1 
           \rmove (-1 0) \fcir f:0 r:0.1 \rmove (-1 0) \fcir f:0 r:0.1 
           \move (7 3) \fcir f:0 r:0.1 
           \rmove (-1 0) \fcir f:0 r:0.1 \rmove (-1 0) \fcir f:0 r:0.1 
           \rmove (-1 0) \fcir f:0 r:0.1 \rmove (-1 0) \fcir f:0 r:0.1 
           \rmove (-1 0) \fcir f:0 r:0.1 \rmove (-1 0) \fcir f:0 r:0.1 
           \rmove (-1 0) \fcir f:0 r:0.1 
           \move (9 2) \fcir f:0 r:0.1 
           \rmove (-1 0) \fcir f:0 r:0.1 \rmove (-1 0) \fcir f:0 r:0.1 
           \rmove (-1 0) \fcir f:0 r:0.1 \rmove (-1 0) \fcir f:0 r:0.1 
           \rmove (-1 0) \fcir f:0 r:0.1 \rmove (-1 0) \fcir f:0 r:0.1 
           \rmove (-1 0) \fcir f:0 r:0.1 \rmove (-1 0) \fcir f:0 r:0.1 
           \rmove (-1 0) \fcir f:0 r:0.1 
           \move (10 1) \fcir f:0 r:0.1 
           \rmove (-1 0) \fcir f:0 r:0.1 \rmove (-1 0) \fcir f:0 r:0.1 
           \rmove (-1 0) \fcir f:0 r:0.1 \rmove (-1 0) \fcir f:0 r:0.1 
           \rmove (-1 0) \fcir f:0 r:0.1 \rmove (-1 0) \fcir f:0 r:0.1 
           \rmove (-1 0) \fcir f:0 r:0.1 \rmove (-1 0) \fcir f:0 r:0.1 
           \rmove (-1 0) \fcir f:0 r:0.1 \rmove (-1 0) \fcir f:0 r:0.1 
           \move (12 0) \fcir f:0 r:0.1 
           \rmove (-1 0) \fcir f:0 r:0.1 \rmove (-1 0) \fcir f:0 r:0.1 
           \rmove (-1 0) \fcir f:0 r:0.1 \rmove (-1 0) \fcir f:0 r:0.1 
           \rmove (-1 0) \fcir f:0 r:0.1 \rmove (-1 0) \fcir f:0 r:0.1 
           \rmove (-1 0) \fcir f:0 r:0.1 \rmove (-1 0) \fcir f:0 r:0.1            
           \rmove (-1 0) \fcir f:0 r:0.1 \rmove (-1 0) \fcir f:0 r:0.1 
           \rmove (-1 0) \fcir f:0 r:0.1 \rmove (-1 0) \fcir f:0 r:0.1 
         \end{texdraw}
         \medskip
         \caption{A Basis of $R/\langle h\rangle +I$.}
         \label{fig:monomials}
       \end{center}
     \end{figure}
     
     For the ``in particular'' part, we  note that
     by Proposition \ref{prop:ordering}
     \begin{displaymath}
       \dim_\C R/\langle h\rangle+I= 
       \dim_\C R/L_{<_{(p,q)}}\big(\langle h\rangle+I\big)\leq
       \dim_\C R/\big\langle y^B\big\rangle+I,
     \end{displaymath}
     and
     the monomials
     $x^\alpha y^\beta$ with $\alpha p+\beta 
     q<pq$ and $\beta<B$ form a $\C$-basis of the
     latter vector space
     (see also Figure \ref{fig:monomials}). Hence,
     \begin{displaymath}
       \dim_\C R/\langle h\rangle +I\leq 
       \sum_{i=0}^{B-1}\big\lceil q-\tfrac{qi}{p}\big\rceil
       =Bq-\sum_{i=1}^{B-1}\lfloor \tfrac{qi}{p}\big\rfloor.
     \end{displaymath}     
   \end{proof}

   \begin{lemma}\label{lem:deg-bound-2}
     Let $g,h\in R$ such that $L_{(p,q)}(g)=x^Ay^B$ and
     $L_{(p,q)}(h)=y^C$, and consider the ideals $J=\big\langle x^Ay^B,y^C, x^\alpha
     y^\beta\;\big|\;\alpha p+\beta q\geq pq\big\rangle$ and $J'=\big\langle g,h, x^\alpha
     y^\beta\;\big|\;\alpha p+\beta q\geq pq\big\rangle$.  Then
     \begin{displaymath}
       \dim_\C R/J'\leq \dim_\C R/J,
     \end{displaymath}
     and if $Ap+Bq\leq pq$
     and $B\leq C\leq p$, then
     \begin{displaymath}
       \dim_\C R/J=
       Ap+Bq-AB-\sum\limits_{i=1}^{A-1}\big\lfloor\tfrac{pi}{q}\big\rfloor
       -\sum\limits_{i=1}^{B-1}\big\lfloor\tfrac{qi}{p}\big\rfloor
       -\sum\limits_{i=C}^{p-1}\min\left\{A,\big\lceil q-\tfrac{Cq}{p}\big\rceil\right\}.
     \end{displaymath}
     Moreover, if $B=0$, then $\dim_\C R/J\leq A\cdot C$.
   \end{lemma}
   \begin{proof}
     By Proposition \ref{prop:ordering}
     \begin{displaymath}
       \dim_\C R/J'\leq
       \dim_\C R/L_{<_{(p,q)}}(J')\leq
       \dim_\C R/J.
     \end{displaymath}
     Let  $I=\big\langle x^\alpha
     y^\beta\;\big|\;\alpha p+\beta q\geq pq\big\rangle$. Then the
     monomials $x^\alpha y^\beta$ with
     $(\alpha,\beta)\in\Lambda=\big\{(\alpha,\beta)\in\N\times\N\;\big|\;\alpha
     p+\beta q< pq\big\}$ form a 
     basis of  $R/I$. 
     Moreover, the monomials $x^\alpha
     y^\beta$ with
     $(\alpha,\beta)\in\Lambda_1\cup\Lambda_2$ are a basis of $J/I$,
     where 
     \begin{displaymath}
       \Lambda_1=\big\{(\alpha,\beta)\in\Lambda\;\big|\;\alpha
     \geq A \mbox{ and }\beta\geq B\big\}
     \end{displaymath}
     and 
     \begin{displaymath}
       \Lambda_2=\big\{(\alpha,\beta)\in\Lambda\setminus\Lambda_1\;\big|\;\beta\geq C\big\}.
     \end{displaymath}
     (See also Figure \ref{fig:monomialbasis}.) This gives rise to the above values for $\dim_\C
     R/J$. 

     \begin{figure}[h]
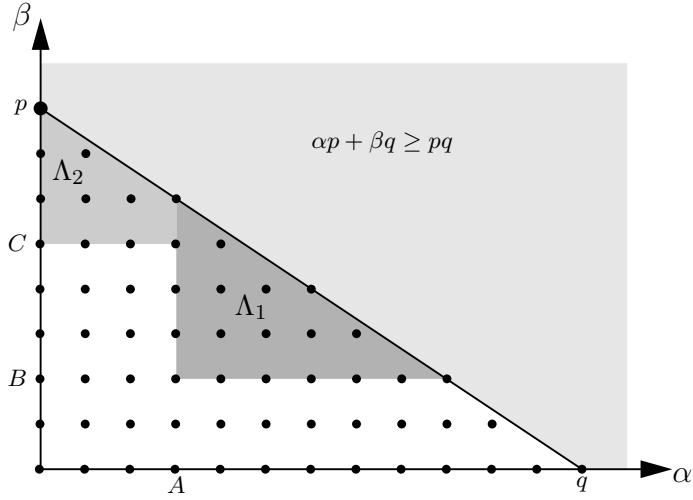

       \bigskip
       \begin{center}
         \begin{texdraw}
           \drawdim cm
           \setunitscale 0.6
           \setgray 0.8
           \arrowheadtype t:F 
           \arrowheadsize l:0.7 w:0.4
           \move (0 8) \lvec (12 0) \lvec (13 0) \lvec (13 9) \lvec (0 9) \ifill f:0.9
           \move (3 2) \lvec (3 6) \lvec (9 2) \ifill f:0.7
           \move (0 5) \lvec (3 5) \lvec (3 6) \lvec (0 8) \ifill f:0.8 
           \setgray 0
           \move (0 0) \avec (0 10) \rmove (-0.2 0)\textref h:R v:C \htext{$\beta$}
           \move (0 0) \avec (14 0) \textref h:L v:T \htext{$\alpha$}
           \move (0 8) \lvec (12 0)
           \move (0 8) \fcir f:0 r:0.15 
           \move (-0.3 8)\textref h:R v:C {\scriptsize\htext{$p$}}
           \move (-0.3 5)\textref h:R v:C {\scriptsize\htext{$C$}}
           \move (-0.3 2)\textref h:R v:C {\scriptsize\htext{$B$}}
           \move (12 -0.2)\textref h:C v:T {\scriptsize\htext{$q$}}
           \move (3 -0.2)\textref h:C v:T {\scriptsize\htext{$A$}}
           \move (4.3 3.3)\textref h:L v:B {\small\htext{$\Lambda_1$}}
           \move (0.25 6.3)\textref h:L v:B {\small\htext{$\Lambda_2$}}
           \move (6 7)\textref h:L v:B {\scriptsize\htext{$\alpha
               p+\beta q \geq pq$}}
           \move (1 7) \fcir f:0 r:0.1 
           \rmove (-1 0) \fcir f:0 r:0.1 
           \move (3 6) \fcir f:0 r:0.1 
           \rmove (-1 0) \fcir f:0 r:0.1 \rmove (-1 0) \fcir f:0 r:0.1
           \rmove (-1 0) \fcir f:0 r:0.1 
           \move (4 5) \fcir f:0 r:0.1 
           \rmove (-1 0) \fcir f:0 r:0.1 \rmove (-1 0) \fcir f:0 r:0.1
           \rmove (-1 0) \fcir f:0 r:0.1 \rmove (-1 0) \fcir f:0 r:0.1
           \move (6 4) \fcir f:0 r:0.1 
           \rmove (-1 0) \fcir f:0 r:0.1 \rmove (-1 0) \fcir f:0 r:0.1 
           \rmove (-1 0) \fcir f:0 r:0.1 \rmove (-1 0) \fcir f:0 r:0.1 
           \rmove (-1 0) \fcir f:0 r:0.1 \rmove (-1 0) \fcir f:0 r:0.1 
           \move (7 3) \fcir f:0 r:0.1 
           \rmove (-1 0) \fcir f:0 r:0.1 \rmove (-1 0) \fcir f:0 r:0.1 
           \rmove (-1 0) \fcir f:0 r:0.1 \rmove (-1 0) \fcir f:0 r:0.1 
           \rmove (-1 0) \fcir f:0 r:0.1 \rmove (-1 0) \fcir f:0 r:0.1 
           \rmove (-1 0) \fcir f:0 r:0.1 
           \move (9 2) \fcir f:0 r:0.1 
           \rmove (-1 0) \fcir f:0 r:0.1 \rmove (-1 0) \fcir f:0 r:0.1 
           \rmove (-1 0) \fcir f:0 r:0.1 \rmove (-1 0) \fcir f:0 r:0.1 
           \rmove (-1 0) \fcir f:0 r:0.1 \rmove (-1 0) \fcir f:0 r:0.1 
           \rmove (-1 0) \fcir f:0 r:0.1 \rmove (-1 0) \fcir f:0 r:0.1 
           \rmove (-1 0) \fcir f:0 r:0.1 
           \move (10 1) \fcir f:0 r:0.1 
           \rmove (-1 0) \fcir f:0 r:0.1 \rmove (-1 0) \fcir f:0 r:0.1 
           \rmove (-1 0) \fcir f:0 r:0.1 \rmove (-1 0) \fcir f:0 r:0.1 
           \rmove (-1 0) \fcir f:0 r:0.1 \rmove (-1 0) \fcir f:0 r:0.1 
           \rmove (-1 0) \fcir f:0 r:0.1 \rmove (-1 0) \fcir f:0 r:0.1 
           \rmove (-1 0) \fcir f:0 r:0.1 \rmove (-1 0) \fcir f:0 r:0.1 
           \move (12 0) \fcir f:0 r:0.1 
           \rmove (-1 0) \fcir f:0 r:0.1 \rmove (-1 0) \fcir f:0 r:0.1 
           \rmove (-1 0) \fcir f:0 r:0.1 \rmove (-1 0) \fcir f:0 r:0.1 
           \rmove (-1 0) \fcir f:0 r:0.1 \rmove (-1 0) \fcir f:0 r:0.1 
           \rmove (-1 0) \fcir f:0 r:0.1 \rmove (-1 0) \fcir f:0 r:0.1            
           \rmove (-1 0) \fcir f:0 r:0.1 \rmove (-1 0) \fcir f:0 r:0.1 
           \rmove (-1 0) \fcir f:0 r:0.1 \rmove (-1 0) \fcir f:0 r:0.1 
         \end{texdraw}
         \medskip
         \caption{A Basis of $R/J$.}
         \label{fig:monomialbasis}
       \end{center}
     \end{figure}

   \end{proof}

   \begin{lemma}\label{lem:degb}
     Let $q>p$ be such that $\frac{q}{p}<\frac{d}{d-1}$ for some
     integer $d\geq 2$, and let
     $0\leq A \leq d$.
     \begin{enumerate}
     \item If $L_{(p,q)}(g)=x^A$,
       then $L_{<_{ds}}(g)=x^A$.
     \item $\m^{p+1}\subseteq \big\langle x^A,y^{p-1},x^\alpha
       y^\beta\;\big|\; \alpha p+\beta q\geq pq\big\rangle$.
     \item If $I$ is an ideal such that $g,h,x^\alpha
       y^\beta\in I$ for $\alpha p+\beta q\geq pq$ and where
       $L_{<_{(p,q)}}(g)=x^A$ and $L_{<_{(p,q)}}(h)=y^{p-1}$, then
       $\degbound(I)\leq p+1$. 
       \\[0.1cm]
       Moreover, if $L_{<_{(p,q)}}(g)$ is minimal among the leading monomials
       of elements in $I$ w.\ r.\ t.\ $<_{(p,q)}$, then $\mult(I)=A$.
     \end{enumerate}
   \end{lemma}
   \begin{proof}
     It suffices to consider the case $A=d$, since this implies the
     other cases. Note that by assumption $d\leq p$. 
     \begin{enumerate}
     \item Since $x^d$ is less than any monomial of degree at least
       $d$ with respect to $<_{ds}$, we have to show that in $g$ no
       monomial of degree less 
       than $d$ can occur with a non-zero coefficient. $x^d$ being the
       leading monomial of $g$ with respect to $<_{(p,q)}$, it suffices to
       show that $\alpha+\beta<d$ implies $\alpha p+\beta q< dp$,
       or alternatively, since $\frac{q}{p}<\frac{d}{d-1}$, 
       \begin{displaymath}
         \alpha + \beta\cdot\frac{d}{d-1}\leq d.
       \end{displaymath} 
       For $\alpha+\beta<d$ the left hand side of this inequality will
       be maximal for $\alpha=0$ and $\beta=d-1$, and thus the
       inequality is satisfied.
     \item 
       We only have to show that $x^\gamma y^{p+1-\gamma}\in
       \big\langle x^d,y^{p-1},x^\alpha y^\beta\;\big|\; \alpha
       p+\beta q\geq pq\big\rangle$ for $\gamma=3,\ldots,d-1$, since the remaining generators of
       $\m^{p+1}$
       definitely are. However, by assumption $\frac{q}{p}<
       \frac{d}{d-1}\leq \frac{\gamma}{\gamma-1}$, and thus
       $\gamma\cdot p+(p+1-\gamma)\cdot q\geq pq$.
     \item By the assumption on $I$ we deduce form (a) and (b) that
       $\degbound\big(L_{<_{ds}}(I)\big)\leq p+1$. However, by Remark
       \ref{rem:degreeordering}
       $\degbound(I)=\degbound\big(L_{<_{ds}}(I)\big)$, which proves
       the first assertion.
       \\[0.1cm]
       Suppose now that $\mult(I)<A$, i.\ e.\ there is an $f\in I$
       such that $\mult(f)\leq A-1$. The considerations for (a) show
       that then $L_{<_{(p,q)}}(f)<x^A$ in contradiction to the
       assumption. 
     \end{enumerate}
   \end{proof}

   \leererpunkt

   \bibliographystyle{amsalpha-tom}
   \bibliography{bibliographie}

\providecommand{\bysame}{\leavevmode\hbox to3em{\hrulefill}\thinspace}
\begin{thebibliography}{Wah74}

\bibitem[Bri77]{Bri77}
Joel Briancon, \emph{Description of {$\Hilb^n{\mathbbm C}\{x,y\}$}}, Inv.\
  math.\ \textbf{41} (1977), 45--89.

\bibitem[BrK86]{BK86}
Egbert Brieskorn and Horst Kn{\"o}rrer, \emph{Plane algebraic curves},
  Birkh{\"a}user, 1986.

\bibitem[GLS97]{GLS97}
{Gert-Martin} Greuel, Christoph Lossen, and Eugenii Shustin, \emph{New
  asymptotics in the geometry of equisingular families of curves}, Internat.\
  Math.\ Res.\ Notices \textbf{13} (1997), 595--611.

\bibitem[GLS00]{GLS00}
{Gert-Martin} Greuel, Christoph Lossen, and Eugenii Shustin, \emph{Castelnuovo
  function, zero-dimensional schemes, and singular plane curves}, J.\ Algebraic
  Geom. \textbf{9} (2000), no.~4, 663--710.

\bibitem[GLS01]{GLS01}
{Gert-Martin} Greuel, Christoph Lossen, and Eugenii Shustin, \emph{The variety
  of plane curves with ordinary singularities is not irreducible}, Intern.\
  Math.\ Res.\ Notes \textbf{11} (2001), 542--550.

\bibitem[GLS05]{GLS05}
{Gert-Martin} Greuel, Christoph Lossen, and Eugenii Shustin, \emph{Singular
  algebraic curves}, Springer, 2005.

\bibitem[GrP02]{GP02}
{Gert-Martin} Greuel and Gerhard Pfister, \emph{A \textsc{Singular}
  introduction to commutative algebra}, Springer, 2002.

\bibitem[Iar77]{Iar77}
Anthony Iarrobino, \emph{Punctual {Hilbert} schemes}, Mem.\ Amer.\ Math.\ Soc.\
  \textbf{10} (1977), no.~188, 1--97.

\bibitem[Kei04]{Kei04}
Thomas Keilen, \emph{Smoothness of equisingular families of curves}, To appear
  in: Trans.\ Amer.\ Math.\ Soc.\ (2004), http:// \!\!www. \!\!mathematik.
  \!\!uni-kl. \!\!de/ \!\!\textasciitilde keilen/ \!\!download/ \!\!Keilen003/
  \!\!Keilen003.ps.gz.

\bibitem[Shu91]{Shu91a}
Eugenii Shustin, \emph{On manifolds of singular algebraic curves}, Selecta
  Math.\ Sov. \textbf{10} (1991), 27--37.

\bibitem[Shu97]{Shu97}
Eugenii Shustin, \emph{Smoothness of equisingular families of plane algebraic
  curves}, Math.\ Res.\ Not. \textbf{2} (1997), 67--82.

\bibitem[Wah74]{Wah74}
Jonathan~M.\ Wahl, \emph{Equisingular deformations of plane algebroid curves},
  Trans.\ Amer.\ Math.\ Soc. \textbf{193} (1974), 143--170.

\end{thebibliography}

\end{document}